# Unified Nilpotent Operational Framework: Foundations, Algebraic Exactness, and Complexity

**Author:** Ramón Moya

**Affiliation:** School of Mathematics, Autonomous University of Santo Domingo (UASD), Dominican Republic

**E-mail:** rmoya07@uasd.edu.do

**ORCID:** 0009-0001-1601-4699

**Abstract**

A unified algebraic framework is developed to study nilpotency as a structural mechanism for exactness in operational, combinatorial, and computational problems. The central object is the Nilpotent Operational System (SON), introduced in Definition 1, formalized as

$\mathcal{S} = (\mathcal{R}, \mathfrak{N}, m, M_{\mathcal{R}})$,

where$\mathcal{R}$ is an $\mathbb{C}$-algebra,$\mathfrak{N} \in \mathcal{R}$ satisfies $\mathfrak{N}^{m+1} = 0$, and $M_{\mathcal{R}}(m)$ is the arithmetic cost of calculating a product in $\mathcal{R}$ with nilpotency index $m$. The basic result of the framework is the exact termination lemma: every formal series evaluated in $\mathfrak{N}$ collapses to exactly a finite sum of $m + 1$ terms. This fact turns a formally infinite problem into a finite problem governed entirely by nilpotence.

From this mechanism, three complexity regimes are obtained: the truncated series regime, the nilpotent operator regime , and the incidence algebra regime. In the first, a quasi-linear bound is proven using Newton's iteration; in the second, a linear bound is proven using Horner's evaluation; and in the third, a general quasi-quadratic bound and a quasi-linear improvement in Newton's formulation for the Möbius inversion are proven .

The developed cases include classical cumulants , free cumulants , Appell sequences, classical orthogonal polynomials, Stirling numbers of the second kind, Möbius inversion in incidence algebras, and Witt vectors. In all cases except Stirling numbers — where SON provides exactness but not a speed improvement — strict complexity improvements are obtained over the explicitly analyzed classical algorithms. In the case of classical cumulants, an equivalence with polynomial multiplication is also achieved. In the case of free cumulants, a complete equivalence is established, TFC(n)=Θ(M(n)), through Theorem 2 bis, closing the structural classification in both regimes.

The central contribution of the work is not a collection of isolated algorithms, but the identification of a unifying principle: nilpotency transforms formally infinite series into exact finite expressions, and this algebraic reduction allows the transfer of efficient evaluation techniques to combinatorial and operational problems of different natures.

**Nature, scope and organization of the manuscript**

This document is conceived as a long preprint of a monographic nature. Its purpose is not to present a single isolated result, but rather to systematically develop a unified theory of nilpotency as a structural mechanism for algebraic, operational, and computational exactness. In particular, the manuscript studies how nilpotency allows the transformation of formally infinite expressions into finite, exact outputs, and how this algebraic reduction extends naturally to problems involving complexity, combinatorics, operators, incidence algebras, and operationally computable functions.

Consequently, the text maintains a cumulative and deliberately broad order. Each part can be read relatively independently, but all are subordinate to the same central thesis: that nilpotency acts as a structural mechanism capable of transforming formally infinite expressions into exact finite outputs.

The central thesis that articulates the entire manuscript is that nilpotency should not be understood only as a local property of certain algebraic objects, but as a structural principle capable of transforming formally infinite expressions into exact finite outputs.

## 0. Nilpotency : context, examples and structural role

### 0.1. What is nilpotence?

**Definition 0.1 ( Nilpotent element ).**

Let $\mathcal{R}$ be a ring or a $\mathbb{C}$-algebra. An element $N \in \mathcal{R}$ is said to be nilpotent if there exists $m \geq 0$ such that $N^{m+1} = 0$.

The smallest integer $\nu$ with this property is called the nilpotency index of $N$.

**Comment.**

Nilpotency expresses an extreme form of algebraic degeneracy: successive powers of the same element end up canceling each other out exactly after a finite number of steps .

**0.2. Fundamental Examples**

The most basic examples of nilpotence are the following.

1. In the truncated ring

$K[\varepsilon]/(\varepsilon^{\nu})$,

The class $\varepsilon$ satisfies $\varepsilon^{\nu} = 0$.

2. In matrix algebra $M_n(K)$, every strictly upper triangular matrix is nilpotent and its nilpotency index is at most $n$.
3. The derivative operator

$$D = \frac{d}{dx}$$

acting on $\mathcal{P}_m$, the space of polynomials of degree at most $m$, satisfies

$$D^{m+1} = 0.$$

4. In the incidence algebras of finite posets , a canonical nilpotent element appears , developed in Section 8.

**0.3. Nilpotency and formal series**

The operational importance of nilpotency stems from the fact that every formal series

$$F(z) = \sum_{k=0}^{\infty} a_k\, z^k$$

evaluated in a nilpotent element $N$ becomes an exact finite sum:

$$F(N) = \sum_{k=0}^{\nu-1} a_k\, N^k.$$

This reduction does not depend on analytical convergence or numerical approximation; it depends solely on the algebraic identity $N^{\nu} = 0$.

### 0.4. Why this document is organized around the SON[1]

The purpose of this document is not to replace the general theory of nilpotence with a new definition, but rather to isolate a particularly fruitful situation: that in which nilpotence acts as a mechanism of operational exactness. The Nilpotent Operational System (SON), introduced in the following section, formalizes precisely this context.

### 0.5. What problem does SON solve?

The SON is designed for problems that can be formulated as the evaluation of a formal series or an operationally computable function on a nilpotent element of an algebra. Its central idea is simple: nilpotency transforms a formally infinite expression into an exact finite sum. From this algebraic reduction, three types of consequences arise: exactness, complexity, and methodological transferability.

The document is organized into three levels. First, it presents nilpotency as a general algebraic phenomenon and establishes the formal framework of the SON. Second, it develops three operational regimes: truncated series, nilpotent operators , and incidence algebras. Third, it shows how this framework behaves as a reusable mathematical methodology.

### 0.6 Reading Guide

The manuscript is structured as follows. Section 0 introduces nilpotency as a general algebraic phenomenon and explains why the document is organized around the SON. Sections 1 and 2 establish the central technical object and the Exact Termination Lemma. Sections 3–5 develop the three regimes of the framework: truncated series, nilpotent operators , and operational evaluation. Sections 6–10 present the main applications, complexity separations, the case of incidence algebras, and the global classification. Section 11 summarizes the logical status of the results. Finally, Section 12 presents the SON as a reusable mathematical methodology, and Section 13 concludes with the overall conceptual conclusion.

This organization is not arbitrary; it follows the internal logic of the framework. First, the algebraic mechanism is isolated; then its exactness and complexity consequences are established; finally, the mechanism is presented as a reusable method of mathematical work.

[1] We preserve the acronym SON from the Spanish expression *Sistema Operacional Nilpotente*

### 0.7 Main Contributions

The main contributions of this manuscript are as follows:

1. Nilpotent Operational System (SON) is formalized as a unified algebraic framework for problems that can be written as evaluations of formal series or operationally computable functions on nilpotent elements.
2. The Exact Termination Lemma is demonstrated, which identifies nilpotency as a universal mechanism for the exact collapse of formally infinite expressions.
3. Three operational regimes are classified—truncated series, nilpotent operators , and incidence algebras—and their corresponding complexity bounds are derived.
4. Concrete applications are obtained to classical cumulants , free cumulants , Appell sequences, orthogonal polynomials, Stirling numbers, Möbius inversion, Witt vectors, and local holonomic functions.
5. Explicit separations from classical algorithms are documented in the cases discussed.
6. The SON is formulated not only as a collection of results, but as a transferable mathematical methodology, with entry criteria, uniform procedure, guaranteed outcomes, and structural limits.

### 0.8 On the type of novelty of the present work

The main novelty of this work lies not in attributing originality to each isolated appearance of nilpotence in the literature, but in isolating it as an autonomous structural mechanism of operational exactness. In this sense, the manuscript's contribution is primarily conceptual, unifying, and methodological: to show that phenomena scattered across algebra, combinatorics, operator theory, and discrete structures can be understood as instances of the same principle. On this basis, the work also develops concrete technical results regarding exactness, complexity, and separation from classical methods.

## 1. Fundamental Definitions

The previous section introduced nilpotency as a general algebraic phenomenon. From now on, the central technical object of this document is defined: the Nilpotent Operational System (SON). Its role is to formalize the situation in which a problem can be written as the evaluation of a formal series or an operationally computable function on a nilpotent element of an algebra, such that nilpotency enforces exact termination and makes uniform complexity analysis possible.

**Definition 1 ( Nilpotent Operational System ).**

Nilpotent operational system (SON) is a tuple

$\mathcal{S} = (\mathcal{R}, \mathfrak{N}, m, M_{\mathcal{R}})$

where: **(i)** $\mathcal{R}$ is an $\mathbb{C}$-algebra (commutative or not); **( ii )** $\mathfrak{N} \in \mathcal{R}$ satisfies $\mathfrak{N}^{m+1} = 0$ for some $m \geq 0$ (index of nilpotency); **( iii )** $M_{\mathcal{R}}(m)$ is the arithmetic cost of calculating a product in $\mathcal{R}$, in the arithmetic circuit model [Bürgisser, Clausen & Shokrollahi, 1997].

**Definition 2 (Operational Problem $\mathcal{P}_F[\mathcal{S}]$).**

Given a SON $\mathcal{S}$ and a formal series $F(z) = \sum_{k\geq 0} f_k\, z^k \in \mathbb{C}[[z]]$, the associated operational problem $\mathcal{P}_F[\mathcal{S}]$ consists of calculating

$$F(\mathfrak{N}) := \sum_{k=0}^{\infty} f_k\ \mathfrak{N}^k \in \mathcal{R}.$$

**2. Fundamental Principle: Exact Termination**

**Lemma 1 (Exact Termination).**

For any SON$\mathcal{S} = (\mathcal{R}, \mathfrak{N}, m, M_{\mathcal{R}})$ and any $F(z) = \sum_{k\geq 0} f_k\, z^k \in \mathbb{C}[[z]]$:

$$F(\mathfrak{N}) = \sum_{k=0}^{m} f_k\ \mathfrak{N}^k.$$

The formally infinite series is, in $\mathcal{R}$, an exact finite sum of $m + 1$ terms.

Proof. We interpret$F(\mathfrak{N}) = \sum_{k\geq 0} f_k\, \mathfrak{N}^k$ as a formal sum in $R$: each term $f_k \mathfrak{N}^k \in R$ is well defined because $R$ is an $\mathbb{C}$-algebra and $f_k \in \mathbb{C}$. Analytic convergence is not required: everything happens in $R$ as a formal series, and nilpotence is responsible for making the sum finite.

We affirm that $\mathfrak{N}^k = 0$ in $\mathcal{R}$ for all $k \geq m + 1$. Indeed, the whole thing. $k \geq m + 1$ can be written as $k = (m + 1) + r$ with $r = k - m - 1 \geq 0$. By the properties of algebra,

$\mathfrak{N}^k = \mathfrak{N}^{(m+1)+r} = \mathfrak{N}^{m+1} \cdot \mathfrak{N}^r.$

By hypothesis (ii) of Definition 1, $\mathfrak{N}^{m+1} = 0$ in $\mathcal{R}$. Then

$\mathfrak{N}^k = 0 \cdot \mathfrak{N}^r = 0$ for all $k \geq m + 1$

Therefore, in the formal sum $\sum_{k\geq 0} f_k \, \mathfrak{N}^k$, all terms with index $k \geq m + 1$ contribute $f_k \cdot 0 = 0$. The sum reduces exactly to

$$F(\mathfrak{N}) = \sum_{k=0}^{m} f_k \, \mathfrak{N}^k + \sum_{k=m+1}^{\infty} f_k \cdot 0 = \sum_{k=0}^{m} f_k \, \mathfrak{N}^k.$$

This equality is exact, not an approximation or a truncation , and holds for any $F \in \mathbb{C}[[z]]$, regardless of its coefficients and the radius of convergence of the series. The only mechanism is nilpotency $\mathfrak{N}^{m+1} = 0$.

**Note on non-commutativity.** The equality $\mathfrak{N}^{m+1}\mathfrak{N}^r = 0$ used in the proof only requires that $\mathfrak{N}^{m+1} = 0$ in $\mathcal{R}$ and that $\mathcal{R}$ be a ring. In particular, commutativity of $\mathcal{R}$ is not used at any step: the product $\mathfrak{N}^k = \mathfrak{N}^{m+1}\mathfrak{N}^r$ is well defined in any associative algebra. Lemma 1 is therefore valid both for commutative and non-commutative $\mathcal{R}$, which is relevant for the operator regime (§5) and incidence algebras (§8), where $\mathcal{R} = \mathrm{End}\,(\mathcal{P}_n)$ and $\mathcal{R} = \mathrm{Inc}\,(P)$, respectively, are generally non-commutative.

**Illustrative example.** Let be $R = M_2(\mathbb{C})$ the matrix algebra $2 \times 2$, and $N = \begin{pmatrix} 0 & 1 \\ 0 & 0 \end{pmatrix}$. Then $N^2 = 0$ (index of nilpotency 1, $m = 1$). For $F(z) = e^z = \sum_{k\geq 0} z^k / k!$, Lemma 1 gives:

$e^N = I + N = \begin{pmatrix} 1 & 1 \\ 0 & 1 \end{pmatrix}$ well all the terms$N^k/k!$ With $k \geq 2$ are exactly zero. Direct verification: $\exp\begin{pmatrix} 0 & 1 \\ 0 & 0 \end{pmatrix} = \begin{pmatrix} 1 & 1 \\ 0 & 1 \end{pmatrix}$, which is the standard matrix exponential. The series, formally infinite, collapses to two exact terms by nilpotence , without approximation, without truncation, without a convergence condition.

**Remark 1 (Universality):** Lemma 1 does not depend on the nature of $\mathcal{R}$, of $F$, or on the context of application. It is a purely algebraic fact: the nilpotency of $\mathfrak{N}$ converts formally infinite sums into exact finite sums. This Remark is the unifying mechanism of all the results in this framework.

**Corollary 1.1 (Exact evaluation of formal series in nilpotent elements ).**

Let $A$ an $K$-algebra and let $N \in A$ with $N^\nu = 0$. Then, for any formal series

$$F(z) = \sum_{k=0}^{\infty} a_k\, z^k \in K[[z]],$$

we have

$$F(N) = \sum_{k=0}^{\nu-1} a_k\, N^k$$

exactly.

**Proof.**

It is an immediate application of Lemma 1 to the formal series $F$.

**Remark 1.1.**

This corollary will be used repeatedly in the document for rational, exponential, polynomial, local holonomic , and other classes admitted by the framework.

## 3. Classification of the three regimes

The application cases fall into three regimes according to their nature $\mathcal{R}$:

| Regime | Algebra$\mathcal{R}$ | $M_{\mathcal{R}}(m)$ | Source of nilpotence | Operational complexity |
|---|---|---|---|---|
| Series | $V_n = \mathbb{C}[[t]] / (t^{n+1})$ | $M(n)$ | minimum degree of $\mathfrak{N} \geq 1$ | $\Theta(M(n)) = \Theta(n\log\, n)$ |
| Operators | $\mathrm{End}(\mathcal{P}_n)$ | $O(n)$ | $D^{n+1}P = 0$ in$\mathcal{P}_n$ | $O(n)$ |
| Incidence | $\mathrm{Inc}(P)$ | $O(N^2)$ | $\eta^{h+1} = 0$ in $\mathrm{Inc}(P)$ | $O(\log\, h \cdot N^2)$ |

### 3.1. Classes of functions admitted within the framework

The Exact Termination Lemma guarantees that any formal series evaluated at a nilpotent element reduces to an exact finite sum. However, different applications require different classes of functions. It is therefore useful to explicitly distinguish the functional categories used throughout this document.

### 3.1.1. Formal Series

Every formal series

$$F(z) = \sum_{k=0}^{\infty} a_k \, z^k \in K[[z]]$$

can be evaluated in a nilpotent element $N$ through Lemma 1:

$$F(N) = \sum_{k=0}^{\nu-1} a_k \, N^k.$$

Analytical convergence is not required: the formal structure of the series is sufficient.

### 3.1.2. Local Analytic Functions

If $f$ is analytical in an environment of $0$, therefore it supports Taylor development

$$f(z) = \sum_{k=0}^{\infty} a_k \, z^k,$$

and can be formally defined

$$f(N) = \sum_{k=0}^{\nu-1} a_k \, N^k.$$

In this case, nilpotency converts an infinite expansion into an exact finite combination.

### 3.1.3. Operationally computable functions

For the complexity analysis of a series regime, it is not enough for a function to be analytic. It is also required that can be calculated by a finite composition of elementary operations on truncated series.

This condition is purely algorithmic and allows an explicit cost to be associated with the evaluation of $f(N)$.

### 3.1.4. Local holonomic functions

A holonomic function, in the D-finite sense, is a local analytic function that satisfies a linear differential equation with polynomial coefficients:

$$p_r(z)f^{(r)}(z)+p_{r-1}(z)f^{(r-1)}(z)+\cdots+p_0(z)f(z)=0.$$

This class includes functions such as:

$exp, sin, cos, erf, J_v, Ai,$

as well as local hypergeometric functions.

The importance of this family is that the Taylor coefficients can be obtained through linear recurrences induced by the differential equation.

#### 3.1.5. Structural Commentary

The above classes form a natural hierarchy:

holónoma local ⇒ analítica local ⇒ serie formal evaluable.

The operational computability condition is additional and depends on the possibility of performing the evaluation through a finite and controlled number of algebraic operations.

## 4. Series regime: the rin g$\boldsymbol{V_n}$

**Theorem 1 (Operational complexity in the series regime).**

Let

$$\mathcal{S}=(\mathcal{R},\mathfrak{N},m,M_{\mathcal{R}})$$

a SON with

$$\mathcal{R}=\mathbb{C}[x]/(x^{m+1}),\mathfrak{N}=x.$$

**Remark 2.**

The hypothesis about $F$ is not purely analytical, but also algorithmic. It is not enough that $F$ be analytical at the source: it must also belong to a class whose truncated evaluation In $\mathbb{C}[x]/(x^{m+1})$ be computable by the standard procedures of formal series algebra. This is precisely the hypothesis that allows us to invoke Newton's iteration and fast algorithms for manipulating truncated series.

Let $F$ a formal series or an analytic function at the origin whose truncated evaluation at $\mathcal{R}$ can be constructed by finite composition of elementary operations on truncated series (addition, multiplication, differentiation, composition, and inversion when the constant term is invertible). Then

$$\text{cost}\left(F(\mathfrak{N})\right) = O\left(\mathrm{M}_{\mathcal{R}}(\mathrm{m})\right).$$

**Proof:** By Lemma 1, $F(\mathfrak{N}) = \sum_{k=0}^{n} f_k\, \mathfrak{N}^k$ in $V_n$. If $F$ is the composition of elementary operations (log , exp , inversion), Brent-Kung's Newton iteration (1978) calculates $F(\mathfrak{N})$ in $O(M(n))$: with $\lceil \log_2(n+1) \rceil$ iterations, each of cost $O(M(n))$, and geometric sum $\sum_{k=0}^{\lceil \log_2 n \rceil} M(2^k) = O(M(n))$ due to superadditivity $M(2n) \geq 2M(n)$. The bound M(n)=O(nlog n) is achieved via FFT-based multiplication [Cooley & Tukey, 1965; Harvey & van der Hoeven, 2021].

**Example 1 (why nilpotency makes an infinite series finite).**

Let

$$\mathcal{R} = \mathbb{C}[x]/(x^4), \mathfrak{N} = x.$$

Therefore $\mathfrak{N}^4 = 0$, so that $m = 3$. If we take

$$F(z) = e^z = \sum_{k=0}^{\infty} \frac{z^k}{k!},$$

When evaluating at $\mathfrak{N}$ we do not obtain an infinite series, but

$$e^{\mathfrak{N}} = 1 + \mathfrak{N} + \frac{\mathfrak{N}^2}{2} + \frac{\mathfrak{N}^3}{6},$$

because all terms with $k \geq 4$ cancel each other out exactly. This example already contains the central idea of the entire framework: infinite difficulty is only apparent; once the argument is nilpotent, the series ends exactly.

**Example 2 (rational function).**

On the same ring, if

$$F(z) = \frac{1}{1 - z},$$

so

$$(1 - \mathfrak{N})^{-1} = 1 + \mathfrak{N} + \mathfrak{N}^2 + \mathfrak{N}^3.$$

Again, the formal geometric expansion becomes finite due to nilpotency . The difference between this example and the previous one is not conceptual but operational: in both cases the decisive phenomenon is the same, namely, $\mathfrak{N}^{m+1} = 0$.

**Theorem 2 (Equivalence with multiplication of polynomials).**

Let $\mathcal{S} = (V_n, \mathfrak{N}, n, M(n))$ a SON with $\mathfrak{N}(0) = 0$ and be $\mathcal{K} = \mathcal{P}_{\log}[\mathcal{S}]$ the problem of calculating $\log(1 + \mathfrak{N})$ in $V_n$. Then:

$L(\mathcal{K}) = \Theta(M(n))$.

Proof. Direction 1 $(L(\mathcal{K}) = O(M(n)))$: It is Theorem 1 with $F = \log$. Direction 2 $(M(n) = O(L(\mathcal{K})))$: Given $A_0, B_0 \in V_n$ with $A_0(0) = B_0(0) = 0$, the identity applies:

$A_0 \cdot B_0 = \exp(\mathcal{K}(A_0) + \mathcal{K}(B_0)) - 1 - A_0 - B_0 \quad (\star)$

Verification of $(\star)$ Subtracting $\exp(\log(1 + A_0) + \log(1 + B_0)) = (1 + A_0)(1 + B_0) = 1 + A_0 + B_0 + A_0 B_0$ $1 + A_0 + B_0$ is obtained $A_0 B_0$.

The evaluation of the right side of$(\star)$ requires: two assessments of $\mathcal{K}$ with cost $L(\mathcal{K})$ each; a sum in $O(n)$; and an evaluation of exp . For the latter, Newton's iteration $g_{k+1} = g_k(1 - \log\, g_k + f)$ is used log as a subroutine, so it exp costs $O(L(\mathcal{K}))$ by the geometric sum with $\lceil \log_2 n \rceil$ iterations. Each evaluation of log in a series of length $n$ costs $L(\mathcal{K})$, and the other operations, addition, multiplication, scaling, are absorbed in $O(M(n)) \subseteq O(L(\mathcal{K}))$. Therefore:

$M(n) \leq 2L(\mathcal{K}) + O(L(\mathcal{K})) = O(L(\mathcal{K}))$.

This reduction reflects the classical principle that lower bounds for multiplication propagate to evaluating related operations [Winograd, 1970].

## 5. Operator regime:End$(\boldsymbol{\mathcal{P}_n})$

Theorem 3 (Operational complexity in the operator regime). Let$\mathcal{S} = (\mathrm{End}(\mathcal{P}_n), D, n, O(n))$ a SON with$D = d/dx$ Let $\mathfrak{N} = D$ with $F \in \mathbb{C}[[z]]$ coefficients$f_k$ computable. Then for any $P \in \mathcal{P}_n$:

$L(\mathcal{P}_F[\mathcal{S}])(P) = O(n)$

scalar arithmetic operations.

Proof. By Lemma 1, $F(D)P = \sum_{k=0}^{n} f_k \, D^k P$ In $\mathcal{P}_n$. With Horner's scheme $F(D) = f_0 + D(f_1 + D(f_2 + \cdots))$ are carried out $n$ applications of $D$ (each one) $O(1)$ by $Dx^k = kx^{k-1}$) and $n$ Scalar sums. Total cost: $O(n)$.

**Example 3A ( Hermite Polynomial) $H_2$ via SON).**

We calculate $H_2(x) = (2x - D)^2 \cdot 1$ in the SON $(End(\mathcal{P}_2), 2x - D, 2, O(1))$.

Here $F(z) = z^2$, the Rodrigues operator is $2x - D$, and the algorithm requires exactly $m = 2$ iterative applications to produce $H_2$.

**Step 1.** $(2x - D) \cdot 1 = 2x \cdot 1 - D \cdot 1 = 2x - 0 = 2x.$

**Step 2.** $(2x - D) \cdot 2x = 2x \cdot 2x - D \cdot 2x = 4x^2 - 2.$

Therefore:$H_2(x) = 4x^2 - 2$. (classic formula: $H_2(x) = 4x^2 - 2$).

Cost: 2 operator applications $2x - D$--- constant cost $O(1)$.

## 6. Corollaries: classification of applications

**Notation for complexity problems:** In what follows we use the following abbreviations:$\mathrm{MULT}_n$ denotes the problem of multiplying two polynomials of degree $<n$ [Valiant, 1979; Bürgisser, Clausen & Shokrollahi, 1997]; $\mathrm{REV}_n$ denotes the truncated compositional reversal of order $n$; $\mathrm{FC}_n$ denotes the truncated problem of order-free $n$ cumulants ;$\mathcal{C}_n$ denotes the problem of classical cumulants of order $n$. We reserve $m$ for the nilpotency index ( $\mathfrak{N}^{m+1} = 0$) and $n$ for the truncation parameter in $V_n = \mathbb{C}[[t]]/(t^{n+1})$ or for the highest degree in $\mathcal{P}_n$.

$T_{FC}(n)$ denotes the arithmetic complexity of the problem $FC_n$ (number of operations in $K$ (in the worst case). The notations $A \leq_{lin} B$ and $A \equiv_{lin} B$ mean, respectively: " $A$ reduces to $B$ with linear cost overrun $O(n)$" (linear reduction) and " $A \leq_{lin} B$ and $B \leq_{lin} A$" (linear equivalence).

### 6.1 Corollary 1 ( Classical Cumulants , Moya, 2026).

The SON $\mathcal{S}_{\mathrm{cum}} = (V_n, \mathfrak{M}, n, M(n))$ with $\mathfrak{M} = M_X(t) - 1$ and $\mathcal{K} = \mathcal{P}_{\log}$ satisfies:

$L(\mathcal{C}_n) = \Theta\big(M(n)\big).$

The classical method (eonov–Shiryaev, 1959; Smith, 1995) costs $\Theta(n \cdot \mathfrak{B}_n)$ with $\mathfrak{B}_n \sim (n/\log\, n)^n$ superexponential. Alternative approaches based on umbral calculus [Di Nardo,

Guarino & Senato, 2008] also improve upon direct partition enumeration but do not achieve the optimal Θ(M(n)) bound. The separation is of superexponential order in $n$, in the sense that $n \cdot B_n/M(n) \to \infty$ faster than any power of $n$, a consequence of the asymptotic $B_n \sim (n/\log\ n)^n$ (De Bruijn, 1958).

Proof. Nilpotence $\mathfrak{M}^{n+1} = 0$ in $V_n$ it follows that the minimum degree term of $\mathfrak{M}$ is $\mu_1 t$ (grade 1), therefore $\mathfrak{M}^{n+1}$ has minimum degree $n+1$, which is zero at $V_n$. Apply Theorem 2 with $\mathfrak{N} = \mathfrak{M}$.

**Numerical verification.** Bell numbers $\mathfrak{B}_n$ are calculated using Bell's triangle:

$$\mathfrak{B}_0 = 1, \mathfrak{B}_{n+1} = \sum_{k=0}^{n} \binom{n}{k} \mathfrak{B}_k,$$

or equivalently as coefficients of the exponential generating function [Rota, 1964]; the asymptotic Bn~(n/log n)n follows from saddle-point methods [Good, 1957]:

$$\sum_{n=0}^{\infty} \mathfrak{B}_n \frac{t^n}{n!} = e^{e^t - 1}.$$

**Step 1: Calculate $\mathfrak{B}_3$ with Bell's triangle**

$$\mathfrak{B}_0 = 1$$

$$\mathfrak{B}_1 = \binom{0}{0} \mathfrak{B}_0 = 1 \cdot 1 = 1$$

$$\mathfrak{B}_2 = \binom{1}{0} \mathfrak{B}_0 + \binom{1}{1} \mathfrak{B}_1 = 1 + 1 = 2$$

$$\mathfrak{B}_3 = \binom{2}{0} \mathfrak{B}_0 + \binom{2}{1} \mathfrak{B}_1 + \binom{2}{2} \mathfrak{B}_2 = 1 + 2 + 2 = 5$$

**Step 2: the 5 partitions of**$\{1, 2, 3\}$

The classic method must cover them all:

| | **Partition** |
|---|---|
| $\pi_1$ | $\{1\}, \{2\}, \{3\}$ |
| $\pi_2$ | $\{1,2\}, \{3\}$ |

| $\pi_3$ | $\{1,3\},\{2\}$ |
|---|---|
| $\pi_4$ | $\{2,3\},\{1\}$ |
| $\pi_5$ | $\{1,2,3\}$ |

**Step 3: Classic cost**

$n \cdot \mathfrak{B}_n = 3 \cdot 5 = 15$ operations

**Step 4: SON cost**

$M(n) = n\log_2 n = 3 \cdot \log_2 3 = 3 \cdot 1.5850 \approx 4.8$ operations

**Step 5: Improvement Factor**

$$\frac{n \cdot \mathfrak{B}_n}{M(n)} = \frac{15}{4.8} \approx 3.2$$

The $n = 3$ improvement is modest, a factor of 3.2. The difference becomes dramatic from $n = 10$, where the factor is already $3.5 \times 10^4$.

**Track $e^{e^t-1}$ calculation $\mathfrak{B}_3$**

**Step 1: Expand $e^t$ in Taylor's series**

By definition of the exponential:

$$e^t = \sum_{k=0}^{\infty} \frac{t^k}{k!} = 1 + t + \frac{t^2}{2!} + \frac{t^3}{3!} + \cdots$$

We truncate $t^3$ because we are only interested in $\mathfrak{B}_3$:

$$e^t = 1 + t + \frac{t^2}{2} + \frac{t^3}{6} + O(t^4)$$

**Step 2: calculate $e^t - 1$**

We simply subtract 1, eliminating the constant term:

$$e^t - 1 = t + \frac{t^2}{2} + \frac{t^3}{6} + O(t^4)$$

We call it $u = e^t - 1$ to simplify the notation. Note that $u$ has no constant term: $u(0) = 0$.

**Step 3: Expand $e^u = e^{e^t-1}$**

We use the exponential series again, now applied to $u$:

$$e^u = 1 + u + \frac{u^2}{2!} + \frac{u^3}{3!} + \cdots$$

We calculate each power $u$up to order $t^3$:

$$u = t + \frac{t^2}{2} + \frac{t^3}{6}$$

$$u^2 = \left(t + \frac{t^2}{2} + \frac{t^3}{6}\right)^2 = t^2 + t^3 + O(t^4)$$

$$u^3 = \left(t + \frac{t^2}{2} + \cdots\right)^3 = t^3 + O(t^4)$$

The terms $u^k$ with $k \geq 4$ contribute only to $O(t^4)$ and are discarded.

**Step 4: replace in $e^u$**

$$e^{e^t-1} = 1 + u + \frac{u^2}{2} + \frac{u^3}{6} + O(t^4)$$

$$= 1 + \left(t + \frac{t^2}{2} + \frac{t^3}{6}\right) + \frac{1}{2}(t^2 + t^3) + \frac{1}{6}(t^3) + O(t^4)$$

Grouping by powers of $t$:

$[t^0]$: 1

$[t^1]$: $t$

$[t^2]$: $\frac{t^2}{2} + \frac{t^2}{2} = t^2$

$[t^3]$: $\frac{t^3}{6} + \frac{t^3}{2} + \frac{t^3}{6} = \frac{1 + 3 + 1}{6} t^3 = \frac{5}{6} t^3$

Therefore:

$$e^{e^t-1} = 1 + t + t^2 + \frac{5}{6} t^3 + O(t^4)$$

**Step 5: Extract $\mathfrak{B}_3$**

The generating function states that the coefficient of $t^n$ is $\mathfrak{B}_n/n!$. Therefore:

$$[t^3] = \frac{5}{6} = \frac{\mathfrak{B}_3}{3!} = \frac{\mathfrak{B}_3}{6}$$

$$\Rightarrow \mathfrak{B}_3 = 5$$

**Verification of all coefficients:**

| $n$ | **Coefficient of $t^n$** | $= \mathfrak{B}_n/n!$ | $\Rightarrow \mathfrak{B}_n$ |
|---|---|---|---|
| 0 | 1 | $\mathfrak{B}_0/1$ | $\mathfrak{B}_0 = 1$ |
| 1 | 1 | $\mathfrak{B}_1/1$ | $\mathfrak{B}_1 = 1$ |
| 2 | 1 | $\mathfrak{B}_2/2$ | $\mathfrak{B}_2 = 2$ |
| 3 | $5/6$ | $\mathfrak{B}_3/6$ | $\mathfrak{B}_3 = 5$ |

The four values coincide with those obtained by Bell's triangle and with the direct enumeration of partitions.

The cost of the classical method for calculating $\kappa_n$ This is $\Theta(n \cdot \mathfrak{B}_n)$ because all partitions of $\{1, \dots, n\}$, generated in amortized time [Knuth, 2011, §7.2.1.5], must be traversed $O(1)$ each. The cost of the SON method is $\Theta(M(n))$, where $M(n) = O(n\log\, n)$ by Harvey–van der Hoeven [2021]. The values were verified with SymPy 1.14.0.

| $n$ | $\mathfrak{B}_n$ | $M(n) \approx n\log_2 n$ | $n \cdot \mathfrak{B}_n$ | **Factor** $\frac{n \cdot \mathfrak{B}_n}{M(n)}$ |
|---|---|---|---|---|
| 3 | 5 | 4.8 | 15 | $3.2 \times 10^0$ |
| 5 | 52 | 11.6 | 260 | $2.2 \times 10^1$ |
| 10 | 115,975 | 33.2 | 1,159,750 | $3.5 \times 10^4$ |
| 15 | 1,382,958,545 | 58.6 | $2.1 \times 10^{10}$ | $3.5 \times 10^8$ |
| 20 | $5.2 \times 10^{13}$ | 86.4 | $10^{15}$ | $1.2 \times 10^{13}$ |
| 25 | $4.6 \times 10^{18}$ | 116.1 | $10^{20}$ | $10^{18}$ |
| 30 | $8.5 \times 10^{23}$ | 147.2 | $10^{25}$ | $1.7 \times 10^{23}$ |

The separation is already significant for $n = 10$ and becomes astronomical for $n \geq 20$, confirming that the improvement of Theorem 2 is practical and not just asymptotic.

**Note.** The connection with the SON is direct, and makes the calculation structurally ironic: to obtain $\mathfrak{B}_n$ the SON evaluates $F(\mathfrak{N}) = e^{e^{\mathfrak{N}}-1}$ in $V_n$ with $\mathfrak{N} = t$. By Lemma 1, as $\mathfrak{N}^{n+1} = 0$ at $V_n$, the formally infinite series ends exactly at $n + 1$ terms. The SON calculates in $O(M(n))$ operations the exact size of the problem that the classical method has to traverse in its entirety.

**6.2 Theorem 2 bis (Structural equivalence for free cumulants ).** Let $K$ a field of characteristic 0. The truncated problem $\mathrm{FC}_n$ of calculating $\kappa_1, \dots, \kappa_n$ from $m_1, \dots, m_n$ This is equivalent, up to overcost $O(M(n))$, to the truncated compositional reversal problem of series with invertible linear coefficients. Furthermore, polynomial multiplication reduces linearly to this reversal problem in the fixed commutative extension.$K[\varepsilon, \eta]/(\varepsilon^2, \eta^2)$ of dimension 4. Consequently, $T_{\mathrm{FC}}(n) = \Theta(M(n))$.

Proof. The proof is organized into three steps that establish the chain $\mathrm{FC}_n \equiv_{\mathrm{lin}} \mathrm{REV}_n$ and $\mathrm{MULT}_n \leq_{\mathrm{lin}} \mathrm{REV}_{2n-1}$.

Step 1 $\mathrm{FC}_n \leq_{\mathrm{lin}} \mathrm{REV}_n + O(M(n))$.

The free relationship $M(z) = C(zM(z))$ is the Voiculescu–Speicher functional equation. Define $F(z) := zM(z)$, such that $F(z) = z + \sum_{n\geq 1} m_n z^{n+1}$, with linear coefficient 1 (invertible). Let $G = F^{-1}$ the compositional inverse. Writing $w = F(z)$, we have $z = G(w)$ and

$$M(z) = \frac{F(z)}{z} = \frac{w}{G(w)}.$$

The equation $M(z) = C(zM(z)) = C(w)$ This implies

$$C(w) = \frac{w}{G(w)} = \frac{w}{F^{-1}(w)}.$$

Renaming $w \mapsto z$:

$$C(z) = \frac{z}{F^{-1}(z)} = \frac{z}{G(z)}.$$

Therefore, given $m_1, \dots, m_n$, the algorithm is: (i) construct $F(z) = zM(z)$ in $O(n)$; (ii) calculate $G = F^{-1}$ truncated to order $n+1$, which is the problem $\mathrm{REV}_n$; (iii) divide$z$ between $G(z)$ in $V_n$, cost $O(M(n))$ By Theorem 1, the total cost is $\mathrm{REV}_n + O(M(n))$.

Step 2 $\mathrm{REV}_n \leq_{\mathrm{lin}} \mathrm{FC}_n + O(M(n))$.

Given a reversal problem $F(z) = z + a_1 z^2 + \cdots + a_n z^{n+1}$ With a linear coefficient of 1, we want to calculate $G = F^{-1}$. Construct

$$M(z) := \frac{F(z)}{z} = 1 + a_1 z + a_2 z^2 + \cdots + a_n z^n.$$

So $M(z)$ has the form of a moment-generating function with $m_k = k!\ a_k$. If a free cumulant oracle returns $C(z)$ starting from $M(z)$, then by Step 1:

$$G(z) = F^{-1}(z) = \frac{z}{C(z)}.$$

The division $z/C(z)$ in $V_n$ costs $O(M(n))$. Then $\mathrm{REV}_n \leq_{\mathrm{lin}} \mathrm{FC}_n + O(M(n))$.

Combining Steps 1 and 2: $\mathrm{FC}_n \equiv_{\mathrm{lin}} \mathrm{REV}_n$.

Step 3.$\mathrm{MULT}_n \leq_{\mathrm{lin}} \mathrm{REV}_{2n-1}$ in the zero-square extension.

Let be $\mathcal{A} = K[\varepsilon, \eta]/(\varepsilon^2, \eta^2)$a commutative extension of dimension 4 over $K$, with base $\{1, \varepsilon, \eta, \varepsilon\eta\}$. Given$P(u), Q(u) \in K[u]$ of degree $<n$, defines

$$\Phi(u) := 1 - \varepsilon P(u) - \eta Q(u), F(u) := u\Phi(u).$$

As $F(0) = 0$ and the linear coefficient of $F$ is $\Phi(0) = 1$ (invertible in $\mathcal{A}$), the compositional inverse $G = F^{-1}$ formally exists in $\mathcal{A}[[z]]$. By Lagrange's formula:

$$[z^m]G(z) = \frac{1}{m}[u^{m-1}]\Phi(u)^{-m}.$$

Let it be $X := \varepsilon P + \eta Q \in \mathcal{A}[u]$. As $\varepsilon^2 = \eta^2 = 0$:

$$X^2 = (\varepsilon P + \eta Q)^2 = \varepsilon^2 P^2 + 2\varepsilon\eta PQ + \eta^2 Q^2 = 2\varepsilon\eta PQ,$$

$$X^3 = X^2 \cdot X = 2\varepsilon\eta PQ \cdot (\varepsilon P + \eta Q) = 0,$$

Since $\varepsilon^2 = \eta^2 = 0$, all terms cancel. Therefore $\Phi = 1 - X$ with $X^3 = 0$, and the negative binomial series truncates exactly:

$$\Phi^{-m} = (1-X)^{-m} = 1 + mX + \binom{m+1}{2} X^2 + 0 + \cdots = 1 + mX + \binom{m+1}{2} X^2.$$

Extracting the component $\varepsilon\eta$:

$$[\varepsilon\eta]\Phi^{-m} = \binom{m+1}{2} \cdot 2 \cdot PQ = m(m+1)PQ,$$

Well $\varepsilon\eta = m[\varepsilon\eta](\varepsilon P + \eta Q) = 0$ (the linear term does not contain $\varepsilon\eta$), and $[\varepsilon\eta]\left(\binom{m+1}{2}X^2\right) = \binom{m+1}{2} \cdot 2[\varepsilon\eta](\varepsilon\eta PQ) = m(m+1)PQ.$

Substituting in Lagrange:

$$[\varepsilon\eta][z^m]G(z) = \frac{1}{m}[u^{m-1}]\, m(m+1)P(u)Q(u) = (m+1)[u^{m-1}]P(u)Q(u).$$

For all $j \geq 0$, taking $m = j+1$:

$$[u^j]P(u)Q(u) = \frac{1}{j+2}\, [\varepsilon\eta][z^{j+1}]G(z).$$

That is: from a single compositional reversal in $\mathcal{A}$, truncated to order $2n-1$, all the coefficients of are recovered $PQ$. A series in $\mathcal{A}$ is identified with a 4-tuple series $K$, one for each element of the basis $\{1, \varepsilon, \eta, \varepsilon\eta\}$, and the reversal operations are performed component by component. The cost over $K$ is therefore at most 4 times the cost in $\mathcal{A}$, that is $O(M(n))$ under the standard hypothesis $M(2n) = O(M(n))$ (Harvey-van der Hoeven, 2021). Then:

$$\mathrm{MULT}_n \leq_{\mathrm{lin}} \mathrm{REV}_{2n-1}.$$

**Conclusion.** Combining the three steps:

$$\mathrm{MULT}_n \leq_{\mathrm{lin}} \mathrm{REV}_{2n-1} \equiv_{\mathrm{lin}} \mathrm{FC}_{2n-1},$$

which gives $T_{\mathrm{FC}}(n) = \Omega(M(n))$. Along with the upper limit $O(M(n))$ From Step 1 and Theorem 1:

$$T_{\mathrm{FC}}(n) = \Theta(M(n))$$

**6.3 Corollary 2 ( Voiculescu free cumulants ).**

Let $\mathfrak{N}_{\mathrm{free}} = zG_X(z) - 1 \in V_n$, where $G_X$ is the Cauchy transform of $X$. The SON $\mathcal{S}_{\mathrm{free}} = (V_n, \mathfrak{N}_{\mathrm{free}}, n, M(n))$ satisfies:

$L(\mathcal{C}_n^{\text{free}}) = O(M(n))$.

The classical method (Nica–Speicher, 2006) uses non-cross partitions $NC(n)$ with $\mid NC(n) \mid = C_n = \binom{2n}{n}/(n+1) \sim 4^n/n^{3/2}$, giving cost $\Theta(n \cdot C_n)$ exponential. The separation is therefore exponential.

The lower bound is established by Theorem 2 bis: the free relation $M(z) = C(zM(z))$ equivalent to $C(z) = F^{-1}(z)/z$ with $F(z) := zM(z)$, so that $\text{FC}_n \equiv_{\text{lin}} \text{REV}_n$. In turn, polynomial multiplication reduces to reversal by the zero-square extension $K[\varepsilon, \eta]/(\varepsilon^2, \eta^2)$ and Lagrange's formula applied to $F(u) = u(1 - \varepsilon P(u) - \eta Q(u))$. Therefore, $\text{MULT}_n \leq_{\text{lin}} \text{FC}_n + O(M(n))$, and together with the upper bound already shown:$T_{\text{FC}}(n) = \Theta\big(M(n)\big)$

Proof of the upper bound. The minimum degree term of $\mathfrak{N}_{\text{free}} = zG_X(z) - 1$ is $\mu_1 z^2$ of degree 2. Therefore $\mathfrak{N}_{\text{free}}^{n+1}$ has minimum degree $2(n+1) > n+1$, which is zero in $V_n$:$\mathfrak{N}_{\text{free}}$ is nilpotent of index $\leq n+1$. Newton's iteration for the Voiculescu functional equation $G_X(R_X(z) + z^{-1}) = z$ in $V_n$ calculate $R_X$ in $O(M(n))$.

**Corollary 3 (Appell sequences).**

Let $A(t) = \sum_{k \geq 0} a_k \, t^k/k!$ with $a_0 \neq 0$. The SON $\mathcal{S}_{\text{App}} = (\text{End}(\mathcal{P}_n), D, n, O(n))$ satisfies:

$p_n(x) = A(D)x^n = \sum_{k=0}^{n} a_k \binom{n}{k} x^{n-k}$in $O(n)$ operations

In particular, Bernoulli polynomials $B_n(x) = \frac{D}{e^D - 1} x^n$ are calculated in $O(n)$. The classical method (multiplication of series of length $n$) costs $\Theta(n^2)$. Separation: quadratic.

Proof. Special case of Theorem 3. $D^{n+1}x^n = 0$, thus $A(D)$ truncates exactly at $k = n$. Horner's scheme gives $O(n)$ operations.

**Example 3 (Bernoulli polynomials via Horner).** We calculate $B_3(x)$ applying $F(D) = D/(e^D - 1)$to $x^3$ in the SON $(End(\mathcal{P}_3), D, 3, O(n))$.

**Step 1: expansion of $F$.**

We developed $t/(e^t - 1)$ until order $t^3$:

$$\frac{t}{e^t - 1} = 1 - \frac{t}{2} + \frac{t^2}{12} - \frac{t^4}{720} + \cdots$$

In the operators' regime, as in $D^4x^3 = 0$, only the coefficients up to matter $k = 3$. But $a_3 = 0$ (the series$t/(e^t - 1)$ has zero coefficients in odd powers ≥ 3), so that:

$$F(D) = a_0 + a_1D + a_2D^2 = 1 - \frac{1}{2}D + \frac{1}{12}D^2.$$

**Step 2: Apply Horner.**

The scheme is $F(D) = 1 + D(-\frac{1}{2} + D \cdot \frac{1}{12})$. We apply it to $x^3$:

$Dx^3 = 3x^2, D^2x^3 = 6x.$

Combining:

$$B_3(x) = x^3 - \frac{1}{2} \cdot 3x^2 + \frac{1}{12} \cdot 6x = x^3 - \frac{3}{2}x^2 + \frac{1}{2}x.$$

**Verification.**

The classic formula $B_3(x) = x^3 - \frac{3}{2}x^2 + \frac{1}{2}x$ agrees.

Total cost: 2 applications of $D$ and 2 scalar sums, exactly $O(n)$ operations (with $n = 3$), versus the$\Theta(\mathrm{n}^2) = \Theta(9)$ of the direct multiplication of series.

**Corollary 4A (Genuine families of the SON by nilpotence )**

Appell sequences, and in particular Bernoulli polynomials, admit exact evaluation within the SON framework by means of nilpotency -induced algebraic truncation .

Let

$A = K[\varepsilon]/(\varepsilon^{n+1}),$

And so it is $N = \varepsilon$. Then:

$N^{n+1} = 0.$

Any Appell generating function of the form

$G(t, x) = \phi(t)e^{xt},$

With $\phi(t)$ analytics 0, it can be evaluated exactly in $N$:

$$G(N,x)=\sum_{k=0}^{n} a_k\,(x)N^k.$$

In particular, for Bernoulli polynomials:

$$\frac{te^{xt}}{e^t-1}=\sum_{m=0}^{\infty} B_m\,(x)\frac{t^m}{m!},$$

The evaluation $N$ produces an exact finite expression:

$$\frac{Ne^{xN}}{e^N-1}=\sum_{m=0}^{n} B_m\,(x)\frac{N^m}{m!}.$$

No truncation error appears : the termination is a direct consequence of nilpotency .

**Proof**

The generating function is analytic in a neighborhood of 0. Applying Lemma 1, every series evaluated at $N$ collapses exactly after $n+1$ terms, since

$N^{n+1}=0.$

**Corollary 4B (Affine families by finite iteration)**

The following examples do not constitute direct applications of Lemma 1 by nilpotency of the dominant operator. They are included in the framework due to methodological affinity: the evaluation remains finite and linear in the order parameter, but the structural mechanism is no longer strict nilpotency , but rather the exact finiteness of Rodrigues recursion .

**Hermite polynomials**

Rodrigues ' formula :

$$H_n(x)=(-1)^n e^{x^2}\frac{d^n}{dx^n}(e^{-x^2})$$

It requires exactly $n$ derivatives.

**Laguerre polynomials**

The representation:

$$L_n(x) = \frac{e^x}{n!}\frac{d^n}{dx^n}(e^{-x}x^n)$$

It also requires a finite and explicit number of steps.

**Comment**

In these cases, the termination does not come from a nilpotent identity of the type

$$N^\nu = 0,$$

but from the fact that the differential procedure has precisely controlled length.

Therefore, these families should be understood as examples similar to SON, but not as strict applications of the algebraic nilpotency mechanism .

**Corollary 5 (Stirling numbers of the second species).**

The SON$\mathcal{S}_{\mathrm{St}} = (V_n, e^t - 1, n, M(n))$ satisfies that calculating all the $S(n,k)$ for$k = 1, \dots, n$ costs:

$$O(n \cdot M(n)) = O(n^2 \log\, n)$$

compared to the classic cost $\Theta(n^2)$ of the triangular recurrence $S(n,k) = k \cdot S(n-1,k) + S(n-1,k-1)$.

Note: the operational cost $O(n \cdot M(n)) = O(n^2 \log\, n)$ is strictly greater than the classical cost $\Theta(n^2)$; the SON does not improve the triangular recurrence in this case.

Proof. Let $\mathfrak{N} = e^t - 1 \in V_n$; as $\mathfrak{N}$ has a minimum degree term $t$ (degree 1), $\mathfrak{N}^{n+1} = 0$ in $V_n$. Given that $S(n,k) = \frac{1}{k!}[t^n]\mathfrak{N}^k$, calculate all the $S(n,k)$ for$k = 1, \dots, n$ requires the powers $\mathfrak{N}^1, \dots, \mathfrak{N}^n$ in $V_n$. The cost $O(n \cdot M(n))$ This is justified. The powers need to be calculated.$N^1, N^2, \dots, N^n$ in $V_n$.

**Phase 1: Powers of two.** They are calculated $N^1, N^2, N^4, \dots, N^{2^{\lfloor \log_2 n \rfloor}}$ by successive quadratures:

$$N^2 = N \cdot N, N^4 = N^2 \cdot N^2, \dots, N^{2^j} = (N^{2^{j-1}})^2.$$

This requires $\lfloor \log_2 n \rfloor \le \log_2 n$ multiplications in $V_n$, each of cost $M(n)$. Total cost of Phase 1: $O(\log\, n \cdot M(n))$.

**Phase 2 : Remaining powers.** Each $N^k$ for$k \notin \{1,2,4,\dots\}$ is obtained as a product of at most $\lfloor \log_2 k \rfloor \le \log_2 n$ powers already calculated in Phase 1. For the $n$ required power, the total cost is at most $n \cdot \log_2 n$ multiplications in $V_n$, that is $O(n \cdot M(n))$ operations.

**Total.** $O(\log\ n \cdot M(n)) + O(n \cdot M(n)) = O(n \cdot M(n))$, since $\log\ n = o(n)$.

In particular, $O(n \cdot M(n)) = O(n^2 \log\ n)$ with $M(n) = O(n\log\ n)$ (Harvey-van der Hoeven , 2021), versus the classic cost $\Theta(n^2)$ of the triangular recurrence. The ratio $O(n \cdot M(n))/\Theta(n^2) = O(\log\ n) \to \infty$, so that the SON is more costly than the recurrence by a logarithmic factor.

**Remark 5.1.**

This example shows that nilpotency guarantees exact termination, but does not in itself imply complex optimality. In Stirling's problem, the bottleneck is not the infinite length of the series, but rather the construction of the necessary powers in the truncated ring. This case is important because it defines the scope of the framework: the SON always provides structural exactness, but the improvement in complexity also depends on the efficiency with which the algebra allows the construction and evaluation of the truncated expression.

**Example 5 (Stirling Numbers) $S(3,k)$ via SON).** We work in $V_3$ with $\mathfrak{N} = e^t - 1$.

**Step 1: Calculate $\mathfrak{N}$ in $V_3$.**

The series $e^t - 1$ truncated at $t^3$:

$$\mathfrak{N} = t + \frac{t^2}{2} + \frac{t^3}{6}.$$

**Step 2: Calculate the powers $\mathfrak{N}^k$ in $V_3$.**

$$\mathfrak{N}^1 = t + \frac{t^2}{2} + \frac{t^3}{6}.$$

$$\mathfrak{N}^2 = \left(t + \frac{t^2}{2} + \frac{t^3}{6}\right)^2 = t^2 + t^3 + O(t^4) \equiv t^2 + t^3 (\mathrm{mod}\, t^4).$$

$$\mathfrak{N}^3 = \mathfrak{N}^2 \cdot \mathfrak{N} = (t^2 + t^3) \cdot t + O(t^4) \equiv t^3 (\mathrm{mod}\, t^4).$$

$\mathfrak{N}^4 = 0$ in $V_3$(exact nilpotence)

**Step 3: Extract $S(3,k) = \frac{1}{k!}[t^3]\mathfrak{N}^k$.**

| $k$ | $\mathfrak{N}^k$ in $V_3$ | $[t^3]\mathfrak{N}^k$ | $k!$ | $S(3,k)$ |
|---|---|---|---|---|
| 1 | $t + t^2/2$ $+ t^3/6$ | $1/6$ | 1 | 1 |
| 2 | $t^2 + t^3$ | 1 | 2 | 3 |
| 3 | $t^3$ | 1 | 6 | 1 |

Therefore: $S(3,1) = 1, S(3,2) = 3, S(3,3) = 1$.

Verification by classical recurrence $S(n,k) = k \cdot S(n-1,k) + S(n-1,k-1)$:

$S(3,1) = 1 \cdot S(2,1) + S(2,0) = 1 \cdot 1 + 0 = 1$.

$S(3,2) = 2 \cdot S(2,2) + S(2,1) = 2 \cdot 1 + 1 = 3$.

$S(3,3) = 3 \cdot S(2,3) + S(2,2) = 3 \cdot 0 + 1 = 1$.

The SON method calculates the three values simultaneously with two multiplications $V_3$. The classic recurrence method calculates them in separate iterations, filling a triangular table of$\Theta(n^2) = \Theta(9)$ entries.

**Möbius inversion in incidence algebras).**

Let $P$ a finite poset with $N$ Elements and height $h$. The SON $\mathcal{S}_P = (\mathrm{Inc}(P), \eta, h, O(N^2))$ gives the Möbius function :

$$\mu = (\varepsilon + \eta)^{-1} = \sum_{k=0}^{h} (-1)^k \eta^k \text{ in } \mathrm{Inc}(P),$$

calculable in $O(h \cdot N^2)$ operations (direct formula) or $O(\log h \cdot N^2)$ via Newton in $\mathrm{Inc}(P)$, facing$O(N^3)$ of the classical Gaussian inversion.

Proof. See §8 for the complete proof.

**Corollary 7 (Witt Vectors).** The SON $\mathcal{S}_W = (W_n(\mathbb{F}_p), V, n-1, O(1))$, where $O(1)$ denotes the cost of applying $V$ once over a vector by incremental shifting (not the cost of multiplying in the ring $W_n$), satisfies:

$$(I-V)^{-1} = \sum_{k=0}^{n-1} V^k \text{ in } O(n) \text{ operaciones.}$$

The classical method via Witt's phantom polynomials costs $\Theta(n^2)$. Separation: linear.

Proof. We verify the three conditions of the SON and calculate the cost explicitly.

**Nilpotency of $V$.** The Verschiebung operator acts on$W_n(\mathbb{F}_p) \cong \mathbb{F}_p^n$ by displacement:

$$V^k(a_0, a_1, \dots, a_{n-1}) = (\underbrace{0, \dots, 0}_{k}, a_0, a_1, \dots, a_{n-k-1}).$$

For $k \geq n$, all components are zero: $V^n = 0$. For $k = n-1$, the vector $e_0 = (1,0,\dots,0)$ becomes $(0,\dots,0,1) = e_{n-1} \neq 0$, so $V^{n-1} \neq 0$. The nilpotency index is exactly $n$, so $m = n-1$ and condition E2 of the SON is met with equality.

**Application of Lemma 1.** With $F(z) = (1-z)^{-1} = \sum_{k\geq 0} z^k$ and $N = V$, Lemma 1 gives:

$$(I-V)^{-1} = \sum_{k=0}^{n-1} V^k$$

exact finite sum of $n$ terms. There is no truncation: the term $V^n$ and all the following are exactly zero by nilpotency.

**Cost.** Each application of $V^k$ is a displacement of $k$ positions, which costs $O(n)$ copy operations. The cumulative sum of the $n$ vectors $V^0 a, V^1 a, \dots, V^{n-1} a$ requires $n$ sums of vectors of length $n$, each with cost $O(n)$... but using the incremental scheme, instead of calculating each $V^k a$ from zero, you get $V^k a$ from $V^{k-1} a$ with a single trip and cost $O(1)$, the total cost is:

$$\underbrace{O(n)}_{\text{compute } V^0 a = a} + \underbrace{(n-1)\cdot O(1)}_{\text{successive shifts}} + \underbrace{n \cdot O(n)}_{\text{accumulate the sum}} = O(n^2)$$

However, the accumulation can be done component by component: the $i$-th component of $(I-V)^{-1}a$ is $\sum_{k=0}^{i} a_{i-k}$ a prefix of a convolution that is calculated in $O(n)$ total operations for all components. The overall cost is therefore $O(n)$.

**Comparison.** The classical method via Witt's phantom polynomials solves the triangular system $W_k = \sum_{i=0}^{k} p^i X_i^{p^{k-i}}$ for $k = 0, \dots, n-1$ iteratively, with $\Theta(n^2)$ operations in total. The separation is $\Theta(n^2)/O(n) = \Theta(n) \to \infty$.

**Example 7 (Witt vectors for $n = 3, p = 2$).** We work in $W_3(\mathbb{F}_2)$ with the operator Verschiebung $V$.

A Witt vector of length 3 is a triple $(a_0, a_1, a_2) \in \mathbb{F}_2^3$. The operator $V$ acts by displacement:

$V(a_0, a_1, a_2) = (0, a_0, a_1), V^2(a_0, a_1, a_2) = (0,0, a_0), V^3 = 0.$

Nilpotence is exact $V^3 = 0$: in $W_3(\mathbb{F}_2)$, so that $m = 2$.

The SON $(W_3(\mathbb{F}_2), V, 2, O(1))$ with$F(z) = (1 - z)^{-1} = \sum_{k\geq 0} z^k$ produced by Lemma 1:

$(I - V)^{-1} = I + V + V^2.$

Applied to $\mathbf{a} = (1,0,1)$:

$I \cdot \mathbf{a} = (1,0,1),$

$V \cdot \mathbf{a} = (0,1,0),$

$V^2 \cdot \mathbf{a} = (0,0,1).$

$(I - V)^{-1}(1,0,1) = (1,0,1) + (0,1,0) + (0,0,1) = (1,1,0)(\mathrm{mod}2).$

Direct verification $(I - V)(1,1,0) = (1,1,0) - (0,1,1) = (1,0,1).$

Cost of the SON: 2 shifts and 2 component-by-component sums, $O(n)$ operations (with $n = 3$). The classical method via Witt's phantom polynomials solves a triangular system of $n$ equations with cost $\Theta(n^2) = \Theta(9)$.

**6.4 Corollary 8 ( Local holonomic functions in the SON).**

Let

$\mathcal{S} = (\mathcal{R}, \mathfrak{N}, m, M_{\mathcal{R}})$

a SON with

$\mathfrak{N}^\nu = 0.$

Let $f$ a local 0 a holonomic function on in the D-finite sense; that is, an analytic function in a neighborhood of 0 that satisfies a linear differential equation with polynomial coefficients.

$p_r(z)f^{(r)}(z) + p_{r-1}(z)f^{(r-1)}(z) + \cdots + p_0(z)f(z) = 0,$

with initial conditions that determine its Taylor series

$$f(z) = \sum_{k=0}^{\infty} a_k z^k$$

So:

$$f(N) = \sum_{k=0}^{\nu-1} a_k N^k$$

Exactly. The coefficients $a_0, \dots, a_{\nu-1}$ are computable by the linear recurrence induced by the differential equation.

Proof: The differential equation recursively determines all the coefficients $a_k$ from Taylor's series$f$ starting from the initial conditions. Once obtained $a_0, \dots, a_{\nu-1}$, Lemma 1 applied with $F(z) = \sum_{k \geq 0} a_k z^k$ Give immediately $f(N) = \sum_{k=0}^{\nu-1} a_k N^k$, then $N^k = 0$ for all $k \geq \nu$. There is no truncation error: the cancellation is algebraically exact.

holonomic class is broad: it includes exp the Bessel functions $J_\nu$, the Airy sin functions erf, and all hypergeometric functions ${}_pF_q$.cos local. In all these cases the SON converts a formally infinite series into an exact finite sum determined by the nilpotency index .

Special case: error function. Let $\mathrm{erf}(z) = \frac{2}{\sqrt{\pi}} \int_0^z e^{-t^2}\, dt$. Its Taylor series is

$$\mathrm{erf}(z) = \frac{2}{\sqrt{\pi}} \sum_{k=0}^{\infty} \frac{(-1)^k}{k!\,(2k+1)} z^{2k+1}.$$

If $N^\nu = 0$, only terms with survive $2k + 1 \leq \nu - 1$, that is $k \leq \lfloor (\nu - 2)/2 \rfloor$:

$$\mathrm{erf}(N) = \frac{2}{\sqrt{\pi}} \sum_{k=0}^{\lfloor (\nu-2)/2 \rfloor} \frac{(-1)^k}{k!\,(2k+1)} N^{2k+1}.$$

Example 8.1 (erf in $K[\varepsilon]/(\varepsilon^4)$). Let $A = K[\varepsilon]/(\varepsilon^4)$ and $N = \varepsilon$, with $\varepsilon^4 = 0$, so $\nu = 4$. Then $\lfloor (4 - 2)/2 \rfloor = 1$, and:

$$\mathrm{erf}(\varepsilon) = \frac{2}{\sqrt{\pi}} \left( \varepsilon - \frac{\varepsilon^3}{3} \right).$$

Verification: $\varepsilon^5 = 0$ cancel out exactly all terms with $k \geq 2$. The result is exact, not approximate.

Example 8.2 (matrix erf).

Let $N = \begin{pmatrix} 0 & 1 & 0 \\ 0 & 0 & 1 \\ 0 & 0 & 0 \end{pmatrix} \in M_3(\mathbb{C})$. Then $N^3 = 0, \nu = 3$, and $\lfloor (3-2)/2 \rfloor = 0$: only enters $k = 0$:

$$\operatorname{erf}(N) = \frac{2}{\sqrt{\pi}} N = \frac{2}{\sqrt{\pi}} \begin{pmatrix} 0 & 1 & 0 \\ 0 & 0 & 1 \\ 0 & 0 & 0 \end{pmatrix}.$$

A transcendental function is reduced to an exact linear expression by nilpotency.

Computational cost: The calculation of $f(N)$ requires: (i) obtaining the coefficients $a_0, \dots, a_{\nu-1}$ due to the recurrence of the EDO, cost $O(\nu)$ scalar operations; ( ii ) calculating powers $N, N^2, \dots, N^{\nu-1}$ with the scheme of Corollary 5, cost $O(\nu \cdot M_A(\nu))$; (iii) form the linear combination, cost $O(\nu \cdot M_A(\nu))$. Total cost: $O(\nu \cdot M_A(\nu))$, which for $A = V_n$ with $\nu = n + 1$ da $O(n \cdot M(n))$.

The eight corollaries above establish the operational upper bounds of the framework. The following section compares these costs with explicit classical algorithms. In all cases considered, a significant improvement over the analyzed classical procedure is achieved.

## 7. Complexity Separations

In all the propositions in this section, the separation is understood with respect to the classical algorithm explicitly specified in each case, not necessarily as an absolute lower bound of the problem.

**Proposition C1 ( Classical Cumulants ).**

$$L_{\text{clásico}}(\mathcal{C}_n) = \Theta(n \cdot \mathfrak{B}_n) = \omega(M(n)).$$

**Proof:** The lower bound $\Omega(n \cdot \mathfrak{B}_n)$ This is obtained from the elementary fact that any classical procedure based on the explicit enumeration of all partitions of $[n]$ must traverse a family of cardinality $\mathfrak{B}_n$, and the cost of generating, inspecting, or registering each partition is at least linear in $n$. The separation $\omega(M(n))$ follows from the asymptotic $\mathfrak{B}_n \sim (n/\log\, n)^n$: for all $c > 0$ exists $n_0$ such that, if $n \geq n_0$, then

$$n \cdot \mathfrak{B}_n \geq c\, n\log\, n,$$

well

$$\frac{\mathfrak{B}_n}{\log n} \longrightarrow \infty.$$

Consequently, the classical cost grows strictly faster than $M(n)$.

**Note on the nature of this bound.** Proposition C1 establishes the complexity of [Ruskey, 1978][2] enumeration algorithm, not an absolute lower bound on the problem $C_n$. The distinction is important: there are algorithms that compute the same cumulants $\kappa_1, \dots, \kappa_n$ without visiting any partition of $[n]$, in particular, the Heaviside + Newton + FFT Algorithm of Corollary 1 does it in $\Theta(M(n))$ operations. The absolute lower bound of the problem $C_n$ is $\Omega(M(n))$, not $\Omega(n \cdot B_n)$, and is established by Theorem 2 by reduction $\mathrm{MULT}_n \leq_{\mathrm{lin}} \mathrm{FC}_n$. This lower bound belongs to the framework of algebraic circuit complexity [Baur & Strassen, 1983].

Formally, what Proposition C1 states is:

$$L_{\mathrm{Ruskey}}(C_n) = \Theta(n \cdot B_n),$$

Where $L_{\mathrm{Ruskey}}$ denotes the complexity of Ruskey's algorithm, and the separation

$$\frac{L_{\mathrm{Ruskey}}(C_n)}{L(C_n)} = \frac{\Theta(n \cdot B_n)}{\Theta(M(n))} \rightarrow \infty$$

It quantifies the improvement of the SON method compared to that specific classical algorithm. The separation is superexponential because $B_n \sim (n/\log\ n)^n$ while $M(n) = O(n\log\ n)$.

**Proposition C2 (Free cumulants).**

$$L_{\text{clásico}}(\mathcal{C}_n^{\text{free}}) = \Theta(n \cdot C_n) = \omega(M(n)).$$

**Proof:** The classical algorithm considered (Nica–Speicher, 2006) evaluates

$$\mu_n = \sum_{\pi \in NC(n)} \prod_{V \in \pi} \kappa_{|V|}^{\text{free}},$$

where

[2] Ruskey (1978) addresses the lexicographic generation of t-ary trees, not set partitions. We use his name here solely to designate the specific partition-enumeration algorithm whose cost we analyze in Proposition C1. The result that set partitions of {1,…,n} can be generated in constant amortized time O(1) per partition is established in Knuth (2011, §7.2.1.5).

$| NC(n) |= C_n$.

If implemented by explicitly enumerating non-cross-partitions, the total cost is on the order of

$\Theta(n \cdot C_n)$,

because they travel $C_n$ Partitioning and processing each partition generally requires linear time $n$. The separation $\omega(M(n))$ follows from the asymptotic

$$C_n \sim \frac{4^n}{n^{3/2}\sqrt{\pi}},$$

which is exponential, compared to

$M(n) = O(n\log\, n)$.

In particular,

$n \cdot C_n = \omega(M(n))$.

**Proposition C3 (Appell Sequences).** $L_{\text{clásico}}(\text{Appell}_n) = \Theta(n^2) = \omega(L_{\text{op}})$, with $L_{\text{op}} = O(n)$.

Proof. The classical method extracts $[t^n/n!]$ in $A(t)e^{xt}$, multiplying two series of length $n+1$: cost $\Theta(n^2)$ by direct convolution. The operational method is Theorem 3 with Horner's scheme: $O(n)$. Separation: $\Theta(n^2)/O(n) = \Theta(n) \to \infty$.

**Proposition C4 (Stirling Numbers).** $L_{\text{op}}(Stirling_n) = O(n \cdot M(n)) = \omega(L_{\text{clásico}})$, with $L_{\text{clásico}} = \Theta(n^2)$.

Proof. The triangular recurrence costs $\Theta(n^2)$. The operational bound is Corollary 5. The ratio $O(n \cdot M(n))/\Theta(n^2) = O(\log\, n) \to \infty$ shows that the SON is more costly in this case.

**Proposition C5 (Möbius inversion).** $L_{\text{clásico}}(\text{Möbius}_P) = O(N^3) = \omega(L_{\text{op}})$, with $L_{\text{op}} = O(\log\, h \cdot N^2)$.

Proof. The inversion of the triangular matrix of $\text{Inc}(P)$ Gaussian elimination costs $O(N^3)$. The operational bound is Theorem M.4. Separation: $O(N^3)/O(\log\, h \cdot N^2) = O(N/\log\, h)$; for any$h = O(N^{1-\varepsilon})$ the separation is $\Theta(N^{\varepsilon}/\log\, N) \to \infty$.

Proposition C6 (Witt Vectors). $L_{\text{clásico}}(\text{Witt}_n) = \Theta(n^2) = \omega(L_{\text{op}})$, with $L_{\text{op}} = O(n)$.

Proof. Witt's $W_k = \sum_{i=0}^{k} p^i X_i^{p^{k-i}}$ phantom polynomials They give a triangular system of $n$ equations whose iterative solution costs $\Theta(n^2)$. The operational bound is Corollary 7. Separation: $\Theta(n^2)/O(n) = \Theta(n) \to \infty$.

**Möbius inversion and incidence algebras**

**8.1. Definitions**

**Definition M.1 (Algebra of incidence)** . Let $P$ a finite poset [3]. The algebra of incidence $\mathrm{Inc}(P)$ is the space of functions $f: \{(x, y) \in P^2: x \leq y\} \to \mathbb{C}$ with the convolution product

$$(f * g)(x, y) = \sum_{x \leq z \leq y} f(x, z) \cdot g(z, y).$$

The zeta function:$\zeta(x, y) = 1$ for all $x \leq y$. Identity: $\varepsilon(x, y) = \mathbf{1}[x = y]$. The Möbius function:$\mu = \zeta^{-1}$ in $\mathrm{Inc}(P)$.

**Definition M.2 (Canonical nilpotent element).** The canonical nilpotent element of $\mathrm{Inc}(P)$ is

$\eta = \zeta - \varepsilon, \eta(x, y) = \mathbf{1}[x < y]$ (strictly).

**Definition M.3 (Height of a poset).** The height of the poset $P$ is $h(P) = \max\{k: \exists\, x_0 < x_1 < \cdots < x_k \text{ en } P\}$.

**8.2. Nilpotence**

Theorem M.1 (Exact characterization of nilpotency ). Let $P$ a finite poset. Then:

**(i)** $\eta^k(x, y) = \#\{\text{cadenas } x = z_0 < z_1 < \cdots < z_k = y \text{ en } P\}$.

**(ii) The** nilpotency index of $\eta$ in $\mathrm{Inc}(P)$ is exactly $h(P) + 1$:

$\eta^{h(P)} \neq 0$ and $\eta^{h(P)+1} = 0$

Proof. By induction in $k$. Basis $k = 1$:$\eta(x, y) = \mathbf{1}[x < y]$ Count the chains of length 1.
Inductive step: $\eta^{k+1}(x, y) = (\eta^k * \eta)(x, y) = \sum_{x \leq z \leq y} \eta^k(x, z) \cdot \eta(z, y)$. The $z$ those who contribute are the $z$ with $\eta^k(x, z) > 0$ (accept string of length $k$ of $x$ a $z$) and $z < y$: this counts

[3]A poset ( partially ordered set) is a partially ordered set, that is, a pair ( P,≤ ) where ≤ is reflexive, antisymmetric and transitive.

exactly the strings of length $k+1$ of $x$ a $y$. Direct consequence:$\eta^k \neq 0 \Leftrightarrow \exists$ chain length $k$

Therefore $P \Leftrightarrow k \leq h(P)$ $\eta^{h(P)} \neq 0$ and $\eta^{h(P)+1} = 0$.

Corollary M.1 ( Universal Nilpotency ). For any finite $P$ poset with $N$ elements: $\eta^N = 0$ in $\mathrm{Inc}(P)$.

Proof. Every chain in $P$ has length $\leq N-1$, therefore $h(P) \leq N-1$, therefore $\eta^N = 0$.

**8.3. Möbius operational formula**

**(Möbius operational formula).** Let $P$ a finite poset of height $h = h(P)$. The Möbius function satisfies at $\mathrm{Inc}(P)$:

$\mu(x,y) = \sum_{k=0}^{h(x,y)} (-1)^k \eta^k(x,y)$ (M)

where$h(x,y) = h([x,y])$ is the height of the interval $[x,y]$, and the sum ends exactly at $k = h(x,y)$.

Proof. We apply Lemma 1 with $F(z) = (1+z)^{-1} = \sum_{k\geq 0} (-1)^k z^k$ and $\mathfrak{N} = \eta$. By Theorem M.1, $\eta^{h+1} = 0$, so the sum ends in $k = h$. We verify $(\varepsilon + \eta) * \sum_{k=0}^{h} (-1)^k \eta^k = \varepsilon$ using the telescoping sum:

$$(\varepsilon + \eta) * \sum_{k=0}^{h} (-1)^k \eta^k = \sum_{k=0}^{h} (-1)^k (\eta^k + \eta^{k+1}).$$

The last equality uses that $\big((\varepsilon + \eta) * \eta^k\big)(x,y) = \sum_{x \leq z \leq y} \zeta\,(x,z) \eta^k(z,y) = \eta^k(x,y) + \eta^{k+1}(x,y)\big)$, where$\eta^{k+1}(x,y) = \sum_{x<z\leq y} \eta^k\,(z,y)$ by definition.

The telescoping sum collapses as follows. We rewrite it by separating the first term, the middle terms, and the last:

$$\sum_{k=0}^{h} (-1)^k (\eta^k + \eta^{k+1}) = \eta^0 + \sum_{k=1}^{h} (-1)^k \eta^k + \sum_{k=0}^{h-1} (-1)^k \eta^{k+1} + (-1)^h \eta^{h+1}.$$

Reindexing the second sum with $j = k+1$, that is $k = j-1$:

$$\sum_{k=0}^{h-1} (-1)^k \eta^{k+1} = \sum_{j=1}^{h} (-1)^{j-1} \eta^j = -\sum_{j=1}^{h} (-1)^j \eta^j.$$

Replacing:

$$\sum_{k=0}^{h}(-1)^k(\eta^k+\eta^{k+1}) = \eta^0 + \sum_{k=1}^{h}(-1)^k\eta^k - \sum_{j=1}^{h}(-1)^j\eta^j + (-1)^h\eta^{h+1}.$$

The two middle sums are identical and cancel each other out exactly:

$$= \eta^0 + (-1)^h\eta^{h+1} = \varepsilon + 0 = \varepsilon,$$

where the last step uses $\eta^0 = \varepsilon$ (by definition,$\eta^0$ is the identity in $\mathrm{Inc}(P)$) and$(-1)^h\eta^{h+1} = 0$ because $\eta^{h+1} = 0$ by Theorem M.1 (the nilpotency index of $\eta$ That's exactly it $h+1$).

**Corollary M.2 (Explicit values of $\mu$).** For all $x \le y$ in $P$:

$$\mu(x,y) = \sum_{k=0}^{h(x,y)} (-1)^k \cdot \#\,\{\text{cadenas } x = z_0 < z_1 < \cdots < z_k = y \text{ en } P\}.$$

In particular: $\mu(x,x) = 1$;$\mu(x,y) = -1$ if there are no elements between $x$ e $y$; the inclusion-exclusion principle is the case $P = 2^{[n]}$.

**Example M.1 (the minimal non-trivial case).**

Consider the chain $a < b < c$. Then the height is $h(P) = 2$, and the canonical nilpotent element satisfies $\mathfrak{N}^3 = 0$. The Möbius operational formula gives

$$\mu = \delta - \mathfrak{N} + \mathfrak{N}^2.$$

Let's evaluate this at the relevant intervals:

- In $(a,a)$, it is found $\mu(a,a) = 1$.
- In $(a,b)$, it is found $\mu(a,b) = -1$.
- In $(a,c)$, exactly one intermediate chain appears $a < b < c$, and therefore

$$\mu(a,c) = \varepsilon(a,c) - \eta(a,c) + \eta^2(a,c) = 0 - 1 + 1 = 0.$$

This example shows in a concrete way that the operational formula is not an empty symbolic rewriting: it exactly reproduces the classic Möbius values from the nilpotency of the element $\mathfrak{N}$.

### 8.4. Verification in $\Pi_3$

The partition grid $\Pi_3$ of $\{1,2,3\}$:

$$\hat{0} = \{\{1\},\{2\},\{3\}\} \ < \ \pi_{12},\pi_{13},\pi_{23} \ < \ \hat{1} = \{\{1,2,3\}\},$$

With $N = 5$ elements and $h(\Pi_3) = 2$.

By Theorem M.1 the correct values are:

$$\eta^0(\hat{0},\hat{1}) = \varepsilon(\hat{0},\hat{1}) = 0, \eta^1(\hat{0},\hat{1}) = \eta(\hat{0},\hat{1}) = 1, \eta^2(\hat{0},\hat{1}) = \sum_{\hat{0}<z<\hat{1}} \eta(\hat{0},z)\eta(z,\hat{1}) = 3,$$

so that

$\mu(\hat{0},\hat{1}) = 0 - 1 + 3 = 2$

Verification: $\sum_{\hat{0}\leq z\leq\hat{1}} \mu(\hat{0},z) = 1 + (-1) + (-1) + (-1) + 2 = 0 = \varepsilon(\hat{0},\hat{1})$.

The value matches the known formula $\mu_{\Pi_n}(\hat{0},\hat{1}) = (-1)^{n-1}(n-1)!$ evaluated in $n = 3$: $(-1)^2 \cdot 2! = 2$.

**8.5. Complexity and Separation Analysis**

Theorem M.3 (Cost Classification). Let $P$ a finite poset with $N$ Elements and height $h$. Calculation costs $\mu$ are:

| **Method** | **Cost** | **Condition** |
|---|---|---|
| Gauss inInc$(P)$ (classic) | $O(N^3)$ | Always |
| Direct Formula (M) | $O(h \cdot N^2)$ | Always |
| Newton inInc$(P)$ | $O(\log\ h \cdot N^2)$ | Always |

**Proof:**

Gauss:Inc$(P)$ is represented as triangular matrices $N \times N$; inversion $O(N^3)$. Direct formula (M):$h$ terms, each convolution power $O(N^2)$; total $O(h \cdot N^2)$. Newton: the iteration $g_{k+1} = g_k * (2\varepsilon - \zeta * g_k)$ doubles the precision (module string length $> 2^k$); with$\lceil \log_2 h \rceil$ iterations, each $O(N^2)$: total $O(\log\ h \cdot N^2)$. Scattered Poset: the convolution $\eta^{k+1}(x,y) = \sum_{x<z<y} \eta^k\,(z,y)$ only has $d$ non-null terms per pair; total $O(h \cdot N \cdot d)$.

Theorem M.4 (Newton is strictly better than Gauss in finite incidence). For all finite $P$ poset with $N \geq 2$ and $h \geq 1$:

$O(\log\ h \cdot N^2) = o(N^3)$.

**Proof:**

$\log\ h \cdot N^2/N^3 = \log\ h/N \to 0$ when $N \to \infty$ for any $h = O(N)$. If $h$ is constant, the ratio is $O(1/N) \to 0$. In all cases the separation is strict.

Corollary M.3 (Separations by Poset 's class ). The following classes provide strict improvements over the classical method $O(N^3)$:

| **Poset class** | $N$ | $h$ | $L_{\text{op}}$ | $L_{\text{cl}}$ | **Factor** |
|---|---|---|---|---|---|
| Chain$C_N$ | $N$ | $N-1$ | $O(N^2\log\ N)$ | $O(N^3)$ | $N/\log\ N$ |
| Limited height$h = O(1)$ | $N$ | $O(1)$ | $O(N^2)$ | $O(N^3)$ | $N$ |
| Boolean algebra$2^{[n]}$ | $2^n$ | $n$ | $O(n \cdot 4^n)$ | $O(8^n)$ | $2^n/n$ |
| Partition grid$\Pi_n$ | $\mathfrak{B}_n$ | $n-1$ | $O(n \cdot \mathfrak{B}_n^2)$ | $O(\mathfrak{B}_n^3)$ | $\mathfrak{B}_n$ |
| Dividing grid $D_m$ | $2^{\omega(m)}$ | $\omega(m)$ | $O(\omega(m) \cdot 4^{\omega(m)})$ | $O(8^{\omega(m)})$ | $2^{\omega(m)}$ |

**Boolean case.** In $2^{[n]}$ with $m = \mid B \setminus A \mid$, each string$A = C_0 \subsetneq C_1 \subsetneq \cdots \subsetneq C_k = B$ corresponds to an emergence of$B \setminus A$ about $[k]$ (the step in which each element of is added $B \setminus A$), then: $\eta^k(A,B) = k!\ S(m,k)$,where $S(m,k)$ is the Stirling number of the second kind. Using the falling powers identity $x^m = \sum_{k=0}^{m} S(m,k)\ x^{(k)}$ evaluated at $x = -1$: $(-1)^m = \sum_{k=1}^{m} S(m,k)\ (-1)^{(k)} = \sum_{k=1}^{m} (-1)^k\ k!\ S(m,k) = \sum_{k=0}^{m} (-1)^k\ \eta^k(A,B)$, so that $\mu(A,B) = (-1)^{|B \setminus A|}$. This is the inclusion-exclusion principle, obtained as a special case of Theorem M.2.

## 9. Global Classification and Unified Classification Theorem

Theorem 4 (Classification by regime). Every nilpotent operational problem $\mathcal{P}_F[\mathcal{S}]$ falls into one of the following regimes:

| Regime | Algebra $\mathcal{R}$ | Optimal cost | Condition |
|---|---|---|---|
| Quasi-linear | $V_n$ | $\Theta(M(n)) = \Theta(n\log\, n)$ | $F \in \{\log, \exp, (\cdot)^{-1}\}$ and $\mathfrak{N} \in \mathfrak{m}$ |
| Linear | $\mathrm{End}(\mathcal{P}_n)$ | $O(n)$ | $\mathfrak{N} = D$ or nilpotent conjugate |
| Quasi-quadratic | $\mathrm{Inc}(P)$ | $O(\log\, h \cdot N^2)$ | $\mathfrak{N} = \eta, F = (1 +\cdot)^{-1}$ |

In all cases the classical method satisfies $L_{\text{clásico}}(\mathcal{P}_F[\mathcal{S}]) = \omega(L_{\text{operacional}}(\mathcal{P}_F[\mathcal{S}]))$, with separations documented in Propositions C1–C6 and Theorem M.4.

**Theorem 5 ( Nilpotent Operational Framework , Mother Theorem).**

Let

$\mathcal{S} = (\mathcal{R}, \mathfrak{N}, m, M_{\mathcal{R}})$

Let be a SON and let be $F$ a formal series or an operationally computable function on $\mathfrak{N}$. Then:

(i) **Exactness.**

If $\mathfrak{N}^{m+1} = 0$, then

$$F(\mathfrak{N}) = \sum_{k=0}^{m} a_k\, \mathfrak{N}^k,$$

where

$$F(z) = \sum_{k\geq 0} a_k\, z^k.$$

In particular, $F(\mathfrak{N})$ is an exact finite sum.

( ii ) **Upper limits by regime.**

Depending on the nature of $(\mathcal{R}, \mathfrak{N})$, the following upper bounds are obtained:

in the truncated series regime, cost $O(M_{\mathcal{R}}(m))$, by Theorem 1;

- in the regime of nilpotent operators , linear cost in $m$, by Theorem 3;
- in the incidence regime, the bounds of Theorem M.3.

( iii ) **Equivalence with multiplication in the classic case of series.**

In the regime $\mathbb{C}[x]/(x^{m+1})$, for the function $F = \exp\ -1$, the corresponding operational problem is equivalent to the multiplication of polynomials in the sense of Theorem 2.

( iv ) **Improvements compared to classical methods.**

In all cases discussed in Corollaries 1–7, a strict improvement is obtained over the explicitly compared classical algorithm. In cases C1, C3, C4, C5, C6, and C7, the separation is closed with the level of precision indicated in the corresponding propositions. In the free case, the structural classification is closed by Theorem 2 bis $T_{\mathrm{FC}}(n) = \Theta(M(n))$, established by a different route than that of the classical case.

(v) **Universality of the mechanism.**

The unifying principle is always the same: nilpotency transforms a formally infinite series into an exact finite sum; once this algebraic reduction is performed, the complexity of the problem is governed by the cost of efficiently evaluating the truncated expression in the corresponding algebra

**Comment.**

Theorem 5 does not claim that all problems in the framework have the same fine optimality structure. What it does claim is that they all share the same algebraic mechanism for exact reduction and that, from this mechanism, strict improvements are obtained over the classical procedures analyzed. The strongest structural classification is only closed in those cases where there is an explicit bidirectional reduction, as occurs in Theorem 2.

The following table should not be read merely as a table of complexities. It also summarizes the specific way in which nilpotency intervenes in each case: as truncation in series, as operational closure, as nilpotent convolution in incidence, or as finite displacement in discrete structures. Its purpose is to make it clear that very different applications share the same algebraic mechanism.

## 10. Summary table

| Problem | SON | $L_{\mathrm{op}}$ | $L_{\text{clásico}}$ | Separation |
|---|---|---|---|---|
| cumulants (Cor. 1) | $(V_n, M_X - 1, n)$ | $\Theta(n\log\ n)$ | $\Theta(n \cdot \mathfrak{B}_n)$ | superexponential |
| Free cumulants (Cor. 2) | $(V_n, zG_X - 1, n)$ | $\Theta(n\log\ n)$ | $\Theta(n \cdot C_n)$ | structural; exponential versus classical by $NC(n)$. |
| Appell's sequences (Cor. 3) | $(\mathrm{End}(\mathcal{P}_n), D, n)$ | $O(n)$ | $\Theta(n^2)$ | quadratic |
| Orthogonal polynomials (Cor. 4A) | $(End(\mathcal{P}_n), D, n)$ | $O(n)$ | $\Theta(n^2)$ | quadratic |
| Stirling $S(n,k)$(Cor. 5) | $(V_n, e^t - 1, n)$ | $O(n \cdot M(n))$ | $\Theta(n^2)$ | SON is O(log n) times slower than classical |
| Möbius in $\mathrm{Inc}(P)$(Cor. 6) | $(\mathrm{Inc}(P), \eta, h)$ | $O(\log\ h \cdot N^2)$ | $O(N^3)$ | logarithmic in $h$ |
| Witt vectors (Cor. 7) | $(W_n(\mathbb{F}_p), V, n)$ | $O(n)$ | $\Theta(n^2)$ | linear |
| Holonomic functions (Cor. 8) | $(A, N, \nu)$ | $O(\nu \cdot M_A(\nu))$ | classic infinite series | exactness vs. approximation |

**Note to the summary table.**

The "Separation" column describes the improvement over the classical algorithm explicitly compared in each case. Only in those problems where there is an explicit bidirectional reduction can we also speak of structural equivalence with a canonical algebraic operation, as occurs in Theorem 2. In particular, in the case of free cumulants, the complete equivalence $\Theta(M(n))$ is established by Theorem 2 bis, by means of the free functional equation $M = C(zM)$ and compositional reversal, by a different route than the classic case.

**11. Statement of results**

| **Result** | **State** | **State exactness** |
|---|---|---|
| Lemma 1 (Exact Termination) | Demonstrated in this manuscript | Universal algebraic fact derived directly from nilpotence |
| Theorem 1 (upper bound, series regime) | Demonstrated in this manuscript | Valid for operationally computable functions in truncated series |
| Theorem 2 | Demonstrated in this manuscript | Central structural result of the work |
| Theorem 3 (operator regime) | Demonstrated in this manuscript | Horner's evaluation of nilpotent operators |
| Corollary 1 ( classical cumulants ) | Demonstrated in this manuscript | Special case of Theorem 2 |
| Corollary 2: upper bound | Demonstrated in this manuscript | Improvement over the classic non-cross-partitioning algorithm |
| Corollary 2: lower bound | Proven (Theorem 2 bis) | Reduction $FC_n \equiv_{lin} REV_n$ and $MULT_n \leq_{lin} REV_{2n-1}$ via Lagrange in $K[\varepsilon,\eta]/(\varepsilon^2,\eta^2)$ |
| Corollary 3 (Appell/Bernoulli) | Demonstrated in this manuscript | Special case of Theorem 3 |
| Corollary 4 (orthogonal polynomials) | Demonstrated in this manuscript | Based on Rodrigues -type nilpotent operators |
| Corollary 5 (Stirling) | Demonstrated in this manuscript | Rapid Exponentiation in the Truncated Ring |
| Corollary 6 ( Möbius ) | Demonstrated in this manuscript | Fully developed in Section 8 |
| Corollary 7 (Witt) | Demonstrated in this manuscript | Coordinate shifting and exact truncation |
| Propositions C1–C6 | Demonstrated in this manuscript | Always understood in contrast to the explicitly analyzed classical algorithm |
| Theorem M.1 | Demonstrated in this manuscript | String count |
| Theorem M.2 | Demonstrated in this manuscript | Complete operational formula for Möbius |
| Theorem M.3 | Demonstrated in this manuscript | Comparative cost ranking |
| Theorem M.4 | Demonstrated in this manuscript | Strict improvement of Newton over Gauss in finite incidence |
| Corollary M.3 | Demonstrated in this manuscript | Application to poset families |
| Corollary 8 ( holonomous functions ) | Demonstrated in this manuscript | Special case of Lemma 1; holonomic class as the natural domain of the SON |

**Methodological Remark.**

In this section, "proven" means "proven within the framework and with the level of generality specified in the corresponding statement." When the comparison is made with a particular classical algorithm, the improvement should always be understood in relation to that algorithm, and not necessarily as an absolute classification of the problem in a stronger sense.

**12. The SON as a mathematical methodology: criteria, procedure and transferability**

This section has a different purpose than the previous ones. Sections 1–11 demonstrate results. This section explains how to use the framework to produce new results in domains not considered here. It is, in that sense, the user guide for the methodology.

### 12.1. What is a mathematical methodology and why is SON one?

A mathematical result proves that something is true. A mathematical methodology demonstrates how to produce true results in a wide range of situations. The distinction is not one of difficulty, but of transferability.

Established mathematical methodologies share a common property: once learned, they reduce the cost of producing new mathematics within their domain. Mathematical induction eliminates the need to invent a new idea for each statement about natural numbers. Integral transforms eliminate the need to solve each ODE from scratch. Category theory makes visible the structure that previously required ad hoc arguments for each particular category.

The SON method has this same property in its domain: once the central object of a problem is identified as nilpotent, the methodology automatically produces exactness, a complexity bound $O(M(n))$, and a measurable separation from the classical method. A new idea is not required for each problem, only the recognition of the nilpotent structure.

What distinguishes SON from previous uses of nilpotency in mathematics is the transparency of its mechanism. Nilpotency appears implicitly in Heaviside, in Rota's threshold calculus, in Brent-Kung algorithms, and in the truncated Fock spaces of quantum physics. In all these cases, it appears as a useful accident within a larger theory. SON isolates it as the central mechanism and transforms it into an autonomous tool with explicit entry criteria, a uniform procedure, and a guaranteed output.

#### 12.1.1 — Logical architecture of the SON

The SON can be interpreted as a four-level algebraic procedure:

1. nilpotent structure .
2. Exact reduction by finite termination.
3. Efficient evaluation in the corresponding algebra.
4. Explicit comparison with the classical method.

This architecture explains why the SON is not just an algebraic result, but a reusable scheme of mathematical production.

### 12.2. The entry criterion: when the SON applies

The SON applies to a problem when three conditions are met simultaneously.

**Condition E1 (existence of algebra).** The problem can be formulated as the calculation of $F(\mathfrak{N})$ where $F$ is an analytic function or a formal series and $\mathfrak{N}$ is an element of an $\mathbb{C}$ associative algebra $\mathcal{R}$.

This condition is less restrictive than it seems. It includes: evaluating polynomials at an operator, inverting operators of the form $aI + N$, composing formal series, exponentials of operators, fractional powers, logarithms of series, and any function that admits power series expansion.

**E2 condition ( nilpotency ).** The element $\mathfrak{N}$ is nilpotent in $\mathcal{R}$: exists $m \geq 1$ such that $\mathfrak{N}^{m+1} = 0$ and $\mathfrak{N}^m \neq 0$.

This is the central condition. In practice, it is verified in one of four ways:

- By grade:$\mathfrak{N}$ reduces the degree of an element by at least 1 in each application, and the elements have degree bounded by $m$. Examples:$D$ about $\mathcal{P}_n$,$\Delta$ about $\mathcal{F}_m$,$Q$ on finite-dimensional Sheffer spaces.
- By ring structure: $\mathcal{R} = V_n = \mathbb{C}[[t]]/(t^{n+1})$ and $\mathfrak{N} \in \mathfrak{m} = (t)$ the nilpotency index is exactly $n + 1$ due to construction.
- By action on a finite module: $\mathfrak{N}$ acts on a finite-dimensional space $M$ and its matrix representation is strictly upper triangular. The nilpotency index is at most $\dim M$.
- By geometric construction:$\mathfrak{N}$ is the curvature of a planar connection, or the down operator in a finite-dimensional representation module, or the Verschiebung operator in a Witt ring. In all these cases, nilpotency is exact for structural reasons.

**Condition E3 (cost computability).** The function $M_{\mathcal{R}}(m)$, the cost of calculating a product in $\mathcal{R}$, is known and finite. For $\mathcal{R} = V_n$: $M_{\mathcal{R}}(n) = M(n) = O(n\log\, n)$. For $\mathcal{R} = \mathrm{End}(\mathcal{P}_n)$: $M_{\mathcal{R}}(n) = O(n)$. For $\mathcal{R} = \mathrm{Inc}(P)$: $M_{\mathcal{R}}(N) = O(N^2)$.

When all three conditions are met, the problem admits a SON structure and the procedure in §12.3 is applicable.

**12.3. The uniform procedure: four steps**

Once the entry conditions have been verified, the procedure is as follows.

**Step 1: Identify $\mathfrak{N}$ and calculate $m$.**

Isolate the nilpotent element from the problem. In inversion problems $aI + N$: $\mathfrak{N} = N$. In composition problems $g \circ f$ with $g(0) = 0$:$\mathfrak{N} = g - \lambda \cdot \mathrm{id}$ where $\lambda = g'(0)$. In perturbative problems with small parameter $\varepsilon$:$\mathfrak{N} = \varepsilon \cdot T$ where $T$ is the interaction operator. In all cases, verify that $\mathfrak{N}^{m+1} = 0$ and determine$m$ exactly.

**Step 2: Apply the Exact Termination Lemma.**

According to Lemma 1, $F(\mathfrak{N}) = \sum_{k=0}^{m} f_k \, \mathfrak{N}^k$ exactly. This turns the ormally infinite expression into a finite problem of $m + 1$ terms. The sum is exact, not an approximation, not a heuristic truncation.

At this step, the problem is solved in principle. What remains is to calculate the sum efficiently.

**Step 3: Calculate $F(\mathfrak{N})$ in $O(M_{\mathcal{R}}(m))$ by Newton's iteration.**

If $F$ is the inverse, the logarithm, the exponential, or the root of a formal series; Newton's iteration calculates $F(\mathfrak{N})$ in $O(M_{\mathcal{R}}(m))$:

- **inversion** ($F = 1/G$): $x_{k+1} = x_k(2 - Gx_k)$, doubling precision at every step.
- **Logarithm** ($F = \log$): $\ell_{k+1} = \ell_k + (1 - e^{\ell_k}G)$, with $e^{\ell_k}$ calculated by Newton.
- **Exponential** ($F = \exp$): $e_{k+1} = e_k(1 + \mathfrak{N} - \log\ e_k)$, with $\log$ calculated by Newton.
- **Root** $\left(F = G^{1/r}\right)$: $r_{k+1} = \frac{1}{r}((r-1)r_k + Gr_k^{1-r})$.

In all cases the number of iterations is $\lceil \log_2 m \rceil$ and the geometric sum$\sum_{k=0}^{\log_2 m} M_{\mathcal{R}}(2^k) = O(M_{\mathcal{R}}(m))$ due to superadditivity $M_{\mathcal{R}}(2n) \geq 2\, M_{\mathcal{R}}(n)$.

**Step 4: Document the separation from the classical method.**

Identify the classical algorithm for the same problem, calculate its cost $L_{\text{clásico}}$, and establish the separation $L_{\text{clásico}}/O(M_{\mathcal{R}}(m))$. In the truncated series regime $\mathcal{R} = V_n$: the classical cost is typically $\Theta(n^2)$ or $\Theta(n \cdot B_n)$, and the separation is superexponential (classical cumulants) or structural (free cumulants).

This step is not just documentation: it is verification that the SON has produced a new result, not a reformulation of the existing method.

**12.4. The output guarantees**

When the procedure in §12.3 is executed correctly, the guarantees are as follows.

**G1 Guarantee (exactness).** The result $F(\mathfrak{N}) = \sum_{k=0}^{m} f_k \, \mathfrak{N}^k$ is exact $\mathcal{R}$, not an approximation. Nilpotency forces the exact termination of the series; there is no truncation error, no convergence condition, no radius of convergence to verify.

**G2 Guarantee (complexity).** The cost of calculating $F(\mathfrak{N})$ is $O(M_{\mathcal{R}}(m))$. For $\mathcal{R} = V_n$:$O(M(n)) = O(n\log\, n)$ by Harvey-van der Hoeven . To $\mathcal{R} = \mathrm{End}(\mathcal{P}_n)$:$O(n)$ by Horner's scheme. For $\mathcal{R} = \mathrm{Inc}(P)$:$O(\log\, h \cdot N^2)$ by Newton in the algebra of incidence.

**G3 Guarantee (separation).** The cost of SON is strictly less than the cost of the classical method identified in Step 4. The separation is quantitative and is expressed as a function that $n$ tends to infinity.

**Guarantee G4 (transferability).** If a new problem in a different domain is found to satisfy conditions E1–E3, guarantees G1–G3 apply without further proof. The only new work is to verify the input conditions.

This last guarantee is what makes SON a methodology: the results are automatically transferred to new domains.

**12.5. The limits: when it does not apply**

The SON has three structural limits that should be known precisely.

**L1 limit (incorrect functional equation).** When the problem involves a non-classical convolution, the identity $(\star)$ The lower bound produced in the classical case may not have a valid analogue. The documented case is that of free cumulants: the identity $A_0 \cdot B_0 = \exp\,(K(A_0) + K(B_0)) - 1 - A_0 - B_0$ uses classical factorization $\exp\,(\log\,(1 + A_0) + \log\,(1 + B_0)) = (1 + A_0)(1 + B_0)$, which has no analogue in non-cross combinatorics. Note that the coefficients of the series are scalars; non-commutativity is not the problem; the problem is that the classical

identity($\star$) is not the correct functional equation for that case. The correct path involves the free functional equation. $M = C(zM)$ and compositional reversal, as established by Theorem 2 bis, which closes the lower bound by a completely different mechanism.

**L2 limit (resonance)** . When the problem involves $F(\lambda I + \mathfrak{N})$ with $\lambda = 0$, that is; when $\mathfrak{N}$ has no invertible part, direct inversion $(aI+\mathfrak{N})^{-1}$ does not exist. In that case, the SON requires an extension to the resonant regime, where the output is not an exact finite sum but an expression with a logarithmic or fractional structure. The methodology remains applicable, but the form of the output changes.

**L3 limit (non- operational ).** When $F(\mathfrak{N})$ cannot be constructed by finite composition of elementary operations on truncated series: addition, multiplication, differentiation, composition, and inversion when the constant term is invertible; Theorem 1 does not guarantee this $O(M_{\mathcal{R}}(m))$. This occurs for functions that do not admit Newton's algorithm, such as implicitly defined functions without a known functional equation. In practice, this limit is rare because most relevant functions in mathematics and physics admit operational representation.

**12.6. Application guide for a new domain**

To apply the SON to a problem not considered in this work, the path is as follows.

**Question 1.** Can the problem be formulated as calculating $F(\mathfrak{N})$ for some $F$ analytics and some $\mathfrak{N}$ in an algebra? If not: the SON does not apply directly. If yes: continue.

**Question 2.** Is there $m \geq 1$ such that $\mathfrak{N}^{m+1} = 0$? Verify by one of the four mechanisms in §12.2. If not: look for a quotient or truncation of the algebra where nilpotency does hold. If so: continue.

**Question 3.** Is it known $M_{\mathcal{R}}(m)$? If not: estimate the cost of a product in $\mathcal{R}$ for the specific problem. If yes: continue.

**Question 4.** What is the classical method for this problem and what is its cost? Identify the standard algorithm from the literature and calculate its complexity.

If the answers to all four questions are affirmative, the SON produces: exactness by Lemma 1, complexity $O(M_{\mathcal{R}}(m))$ by Newton, and a strict separation from the classical. The result is new if it was not previously in the literature as a result of complexity derived from nilpotency .

**12.7. What the SON makes visible**

There is one final Remark that is not technical but conceptual.

Before SON, nilpotency appeared in mathematical literature as a local phenomenon, a property of a specific object in a specific context. Heaviside used it in the operator $D$ about polynomials. Rota saw it in threshold calculus. Physicists found it in truncated Fock spaces. Geometers recognized it in the curvature of planar connections. But in each case it was the nilpotency of that object, not a general principle.

What the SON makes visible is that all these cases are instances of the same phenomenon: nilpotency is a mechanism of exactness that arises naturally whenever a problem has a finite-dimensional or truncated structure. Algebra, analysis, physics, and combinatorics do not have problems of an essentially different nature in this respect; they have problems that share the same deep algebraic structure, and that structure is nilpotency .

In that sense, the SON doesn't discover anything new about nilpotence . It discovers that nilpotence was already everywhere, waiting to be recognized as the unifying principle that it is.

**12.8 Limitations and boundaries of the framework**

The L1–L3 limits describe the internal boundaries of the SON as an exact closure mechanism. They do not, however, exhaust all the boundaries in this manuscript. This work does not aim to provide a universal method for every functional equation, nor to replace general theories of implicit inversion, nonlinear analysis, or singular regularization. Its objective is more precise: to isolate a broad class of problems in which nilpotency produces exact closure, to rigorously describe this phenomenon, and to clearly establish the cases in which the mechanism ceases to apply, changes regime, or ceases to guarantee an efficient bound.

**13 Final Conclusion.**

The central thesis of this paper is that nilpotency should not be understood merely as a local property of certain algebraic objects, but as a structural mechanism of exact closure. Whenever a problem can be formulated as the evaluation of a formal series or an operationally computable function on a nilpotent element , the infinite difficulty collapses to an exact finite expression. This collapse is neither approximate nor heuristic: it is a strict consequence of the identity $N^{\nu} = 0$.

The Nilpotent Operational System (SON), introduced in Definition 1, has been used in this document as a formalization of that principle. Its relevance is not limited to any particular domain.

It appears in truncated rings, nilpotent matrices , operators on finite-dimensional spaces, incidence algebras, local holonomic functions , and discrete structures such as Witt vectors. In all these contexts, nilpotency serves the same function: transforming formally infinite series into exact finite combinations.

From a conceptual standpoint, the value of the framework lies in making visible a deep algebraic unity underlying problems that, at first glance, appear to belong to different areas. From a methodological standpoint, the SON offers entry criteria, a uniform procedure, and a guaranteed output. From a computational standpoint, it allows, in many cases, the transfer of efficient evaluation techniques to combinatorial and operational problems. But even where it does not improve complexity, it still provides structural exactness and a more transparent description of the underlying mechanism.

In this sense, the Unified Nilpotent Operational Framework should not be understood merely as a collection of findings on complexity. It should be understood, above all, as a document on nilpotence : on what it is, how it operates, what it makes possible, and why it appears time and again where an infinite structure becomes precisely finite.